\documentclass[11pt]{article}
\usepackage{array}
\usepackage{amsmath,amsthm,amssymb,bm}
\usepackage{enumerate, enumitem}
\usepackage{hyperref}
\usepackage{tikz}
\usepackage[mathlines]{lineno}
\usepackage{soul} 
\usepackage{color}
\usepackage{mathrsfs}
\usepackage{mathdots}
\usepackage{bbm}
\usepackage{dsfont}
\usepackage{float}
\usepackage{nicematrix}

\theoremstyle{plain}
\newtheorem{theorem}{Theorem}[section]
\newtheorem{lemma}[theorem]{Lemma}
\newtheorem{proposition}[theorem]{Proposition}
\newtheorem{corollary}[theorem]{Corollary}
\newtheorem{observation}[theorem]{Observation}

\theoremstyle{definition}
\newtheorem{example}[theorem]{Example}
\newtheorem{remark}[theorem]{Remark}

\newtheorem{question}[theorem]{Question}
\newtheorem{problem}[theorem]{Problem}

\DeclareMathOperator{\tr}{tr}
\DeclareMathOperator{\diag}{diag}
\DeclareMathOperator{\adj}{adj}
\DeclareMathOperator{\col}{col}
\DeclareMathOperator{\row}{row}
\DeclareMathOperator{\supp}{supp}

\definecolor{red}{rgb}{.8,0,0}

\definecolor{bblu}{rgb}{0,0,1}

\definecolor{gray}{rgb}{.4,.4,.4}

\definecolor{gre}{rgb}{0,.7,0}

\definecolor{olive}{rgb}{.502,.502,0}

\newcommand{\eref}[1]{(\ref{#1})}

\newcommand{\R}{{\mathbb R}}
\newcommand{\C}{{\mathbb C}}

\newcommand{\U}{\mathcal U}

\newcommand{\ii}{\mathbbm{i}}

\newcommand{\bu}{\mathbf{u}}
\newcommand{\bv}{\mathbf{v}}

\newcommand{\be}{\mathbf{e}}
\newcommand{\bzero}{\boldsymbol{0}}
\newcommand{\bc}{\mathbf{c}}
\newcommand{\br}{\mathbf{r}}
\newcommand{\bw}{\mathbf{w}}

\newcommand{\ac}{\operatorname{K}}
\newcommand{\spn}{\operatorname{span}}
\renewcommand{\Re}{\operatorname{Re}}
\renewcommand{\Im}{\operatorname{Im}}
\newcommand{\parens}[1]{\left( #1\right)}
\newcommand{\rank}{\operatorname{rank}}
\newcommand{\spec}{\operatorname{spec}}
\newcommand{\set}[1]{\left\{#1\right\}}

\allowdisplaybreaks

\begin{document}

\title{Matrix Apportionment: General Cases and Apportionment Constants}
\author{Dustin R. Baker\thanks{Department of Mathematics, Iowa State University, Ames, IA 50011, USA (bakerd@iastate.edu, jmiller0@iastate.edu, pungello0@iastate.edu)} \and Bryan A.~Curtis\thanks{Air Force Research Laboratory, Wright-Patterson Air Force Base, OH 45433, USA (dr.bryan.curtis.math@gmail.com)} \and Joe Miller\footnotemark[1] \and  Hope  Pungello\footnotemark[1]}

\maketitle\vspace{-15pt}

\begin{abstract} 
A matrix is apportionable if it is similar to a matrix whose entries have equal moduli. This paper shows that all nilpotent matrices and all matrices with rank at most half their order are apportionable. General results are established and applied to classify all apportionable matrices of order 2 and partially those of order 3. Additionally, the study of the set of apportionment constants for matrices is initiated.
\end{abstract}



\section{Introduction}\label{section: intro}

The concept of matrix apportionment was recently introduced and explored by Clark et al.\ in \cite{CCGH2024}. In many ways, matrix apportionment is antithetical to Jordan decomposition. While the latter seeks to concentrate the nonzero entries of a matrix near the diagonal, apportionment aims to transform a matrix into one whose entries have equal moduli. More precisely, a square matrix $A$ is \emph{apportionable} if there exists a nonsingular matrix $M$ such that $MAM^{-1}$ is \emph{uniform}, meaning all entries of $MAM^{-1}$ share the same modulus. The results in \cite{CCGH2024} primarily focused on matrices that are unitarily similar to a uniform matrix, termed \emph{$\U$-apportionable}. The present paper delves into the more general case in greater detail.

Throughout this paper, unless otherwise specified, all $m\times n$ matrices $A$ are assumed to have complex entries, i.e., $A\in\C^{m\times n}$. For a square matrix $A$ of order $n$ and a nonsingular matrix $M$ of the same order, if $MAM^{-1}$ is uniform, the common modulus $\kappa$ of its entries is called an \emph{apportionment constant} of $A$. We denote the set of all apportionment constants of $A$ by $\ac(A)$. While some matrices, such as the zero matrix, may have a unique apportionment constant, this is often not the case for others, as is illustrated in many of the following results. When $M$ is unitary we say that $\kappa$ is the \emph{unitary apportionment constant} of $A$. As noted in \cite{CCGH2024}, unitary apportionment constants are unique because the entries of an $n\times n$ uniform matrix $B$ all have the common modulus $\frac{\|B\|_F}{n}$, and the Frobenius norm $\|B\|_F$ is unitarily invariant. 

Given the recent emergence of apportionable matrices, their applications in the existing literature are still relatively few. One significant connection to uniform matrices lies with complex Hadamard matrices. A \emph{complex Hadamard matrix} is defined as a unitary matrix where every entry has a modulus of $1$. Section~\ref{sec: tools} of this paper utilizes the well-known Hadamard inequality to derive bounds on the apportionment constants of apportionable matrices. For further details on Hadamard matrices, the reader may consult \cite{BEHC2021, TZ2006}. An intriguing application of $\U$-apportionable matrices from \cite{CCGH2024} involves their use in studying graceful labelings of graphs. The authors of \cite{CCGH2024} also discuss connections between equiangular lines and uniform matrices. Other research pertaining to equiangular lines has also made use of uniform matrices. Notably, an equiangular tight frame (ETF) is said to be \emph{flat} if every entry has modulus 1 (see, for example, \cite{FJMP2021, GR2009}). Flat ETFs have connections to coding theory and waveform design.

Before presenting the main results, Section~\ref{subsec: background and notation} provides a review of our notational conventions and definitions. 
Section~\ref{sec: tools} establishes preliminary results and addresses \cite[Question 7.3]{CCGH2024}: If $B$ is uniform, is $\spec(B)\cup\{0\}$ the spectrum of a uniform matrix? In Section~\ref{sec: nilpotent}, we demonstrate that all nilpotent matrices $A$ are apportionable and their set of apportionment constants is $\ac(A) = (0,\infty)$. Section~\ref{sec: half rank} proves that all $n\times n$ matrices with rank at most $\frac n2$ are apportionable. Section~\ref{sec: small order} presents additional general results, with a primary focus on applying these findings to classify apportionable matrices of small order. We conclude with Section~\ref{sec: conclusion} where we provide direction for future work.

\subsection{Notational Conventions and Definitions}\label{subsec: background and notation}
We write $\R$ and $\C$ to denote the set of real and complex numbers, respectively. To avoid confusion with indices, we write $\ii$ for the imaginary unit. We will use $\exp(x)$ and $e^x$ interchangeably. The conjugate and modulus of $z\in\C$ are $\overline{z}$ and $|z|$, respectively.

The $n$-dimensional complex vector space is written as $\C^n$ and the set of $m\times n$ complex matrices is denoted $\C^{m\times n}$. All vectors are written in bold Roman font, e.g., $\bv$. The transpose of a vector $\bv$ and matrix $A$ is written $\bv^\top$ and $A^\top$, respectively. The standard basis vectors of $\C^n$ are $\be_1,\ldots,\be_n$ and the zero vector is $\bzero$. We imbue $\C^n$ with the standard dot product defined by $\bv\cdot \bw = \bw^*\bv$. The orthogonal complement of a subspace $V$ of $\C^n$ is denoted $V^\perp$.

The zero matrix is denoted $O$, with its dimension specified as $O_m$ in the square case or $O_{m\times n}$ in the rectangular case. The identity matrix is denoted $I$ with its order specified as $I_m$ as necessary. An $n\times n$ diagonal matrix $D = [d_{i,j}]$ may be written as $D = \diag(d_{1,1},\ldots,d_{n,n})$. We write $E_{i,j}$ for the matrix with a $1$ in the $(i,j)$-entry and $0$ elsewhere; its dimension will be clear from context. The conjugate transpose of a matrix $A$ is denoted $A^*$. The column space and row space of a matrix $A$ are denoted $\col(A)$ and $\row(A)$, respectively. The adjugate of a nonsingular matrix $M$ is $\adj(M) = \det(M) M^{-1}$.

The spectrum of a matrix $A$, denoted $\spec(A)$, is the multiset of its eigenvalues. The spectral radius of $A$ is $\rho(A) = \max\{|\lambda| : \lambda\in\spec(A)\}$. We write $\rank(A)$ to denote the rank of a matrix $A$ and if $\rank(A) = r$ we say $A$ is a rank $r$ matrix. The Jordan form of a matrix $A$ is denoted $J(A)$ and the $k\times k$ Jordan block corresponding to an eigenvalue $\lambda$ is denoted $J_k(\lambda)$.

We now recall the definitions in the introduction related to apportionment. Let $A\in\C^{n\times n}$ and write $A = [a_{i,j}]$. Then $A$ is \emph{uniform} if there exists a nonnegative real number $\kappa$ such that $|a_{i,j}| = \kappa$ for all $i,j\in\{1,\ldots,n\}$. Note that this definition clearly extends to the rectangular case. The matrix $A$ is \emph{apportionable} if there exists a nonsingular matrix $M$ such that $MAM^{-1}$ is uniform, and is $\U$\emph{-apportionable} if there exists a unitary matrix $U$ such that $UAU^*$ is uniform. If $MAM^{-1}$ is uniform, the common modulus $\kappa$ of its entries is called an \emph{apportionment constant} of $A$. The set of all apportionment constants of $A$ is denoted $\ac(A)$.

\section{Preliminary Results}\label{sec: tools}

This section contains tools that are broadly useful in the study of apportionable matrices. One elementary observation is that we may always assume, without loss of generality, that a matrix $A$ is in Jordan form since $A$ is similar to $J(A)$. Explicitly finding a nonsingular matrix $M$ that apportions a matrix $A$ is often made much easier by utilizing the sparsity of $J(A)$. Note that if $\spec(A) = \{\lambda_1,\ldots,\lambda_n\}$, then the Jordan form of $A$ can be written as 
\begin{equation}\label{eqn: Jordan outer products}
J(A) = \sum_{k=1}^n \lambda_k\be_k\be_k^\top + \sum_{k=1}^{n-1} \alpha_k \be_k \be_{k+1}^\top,
\end{equation}
where each $\alpha_k\in\{0,1\}$. The following observation is elementary but provides an explicit way to simplify the Jordan form of matrices for apportionment.

\begin{observation}\label{obs: Jordan form}
Let $A\in \mathbb{C}^{n\times n}$ be a matrix in Jordan form and $\lambda\neq 0$. The Jordan form of $\lambda A$ can be obtained by scaling the diagonal entries of $A$ by $\lambda$. In particular, $S(\lambda A)S^{-1}$ is in Jordan form, where $S=\diag(1,\lambda,\ldots, \lambda^{n-1})$.
\end{observation}

Our first result, Lemma~\ref{lem:pad-by-0}, establishes that if $A$ is apportionable, then $A\oplus O$ is apportionable. We remark that Lemma~\ref{lem:pad-by-0} also answers \cite[Question 7.3]{CCGH2024} to the affirmative, i.e., if $B$ is uniform, then $\spec(B)\cup\{0\}$ is the spectrum of a uniform matrix.

\begin{lemma}\label{lem:pad-by-0} 
Let $A\in\C^{n\times n}$. If $A$ is apportionable, then $A \oplus O$ is apportionable. Moreover, if $\kappa$ is an apportionment constant of $A$, then $\kappa$ is an apportionment constant of $A\oplus O$.
\end{lemma}

\begin{proof}
Assume $A$ is apportionable. It suffices to prove $A\oplus\left[\, 0\, \right]$ is apportionable. Let $M$ be a nonsingular matrix such that $B = MAM^{-1}$ is uniform with common modulus $\kappa\in\ac(A)$. Let $\omega = e^{\pi\ii/3}$ and $D = \diag(\overline{\omega},1,\ldots,1)$. Consider the $(n+1)\times (n+1)$ matrix
\[
N = \left[\begin{array}{c|c}
I & -\omega \be_1 \\
\hline\\[-1em]
\omega \be_1^\top & 1
\end{array}\right]
\quad \text{with inverse} \quad 
N^{-1} = \left[\begin{array}{c|c}
D & \be_1 \\
\hline\\[-1em]
-\be_1^\top & \overline{\omega}
\end{array}\right].
\]
Let $\bc_1$ and $\br_1^\top$ denote the first column and the first row of $B$, respectively. Observe that the matrix
\[
N\left[\begin{array}{c|c}
M & \bzero \\
\hline\\[-1em]
\bzero^\top & 1
\end{array}\right]
\left[\begin{array}{c|c}
A & \bzero \\
\hline\\[-1em]
\bzero^\top & 0
\end{array}\right]
\left( N\left[\begin{array}{c|c}
M & \bzero \\
\hline\\[-1em]
\bzero^\top & 1
\end{array}\right]\right)^{-1}
= 
\left[\begin{array}{c|c}
BD & \bc_1 \\
\hline\\[-1em]
\omega \br_1^\top D & \omega \br_1^\top\be_1
\end{array}\right]
\]
is uniform with common modulus $\kappa$. Thus, $A\oplus\left[\, 0\, \right]$ is apportionable with apportionment constant $\kappa$. 
\end{proof}

To motivate the next lemma, consider $n\times n$ matrices $M_1$ and $M_2$ with columns $\bu_1,\ldots,\bu_n$ and rows $\bv_1^\top,\ldots,\bv_n^\top$, respectively. If $A\in\C^{n\times n}$ is in Jordan form and $\spec(A) = \{\lambda_1,\ldots,\lambda_n\}$, then by \eref{eqn: Jordan outer products}
\[
M_1AM_2 
= M_1\left(\sum_{k=1}^{n} \lambda_k \be_{k} \be_{k}^\top + \sum_{k=1}^{n-1} \alpha_k \be_{k}\be_{k+1}^\top\right)M_2 \\
= \sum_{k=1}^{n} \lambda_k \bu_{k}\bv_{k}^\top + \sum_{k=1}^{n-1} \alpha_k \bu_{k}\bv_{k+1}^\top,
\]
where each $\alpha_k\in\{0,1\}$. Thus, if $0$ is an eigenvalue of $A$, some of the columns of a nonsingular matrix $M$ and some of the rows of $M^{-1}$ may be inconsequential in determining if $MAM^{-1}$ is uniform. In order to leverage this, we make use of the next result which only requires basic linear algebra knowledge.  As a formal reference could not be found, a proof is provided below.

\begin{lemma}\label{lem:block-matrix-inverse}
Let $U \in \C^{n \times m}$ and $V \in \C^{m \times n}$ with $m < n$ and assume  $VU = I_m$. Then there exist matrices $U'\in\C^{n\times(n-m)}$ and $V'\in\C^{(n-m)\times n}$ such that 
$\left[\begin{array}{c} 
V \\ 
\hline 
V' 
\end{array}\right]
\left[\begin{array}{c|c} 
U & U' 
\end{array}\right] 
= I_n$.
Moreover, the columns of $U'$ may be chosen to be any basis of $\col(V^*)^\perp$.
\end{lemma}
\begin{proof}
Let $U'$ be a matrix whose columns form a basis for $\col(V^*)^\perp$. Note that, for any matrix $A$ with suitable dimensions, the $(j,k)$-entry of $VA$ is the inner product between column $k$ of $A$ and column $j$ of $V^*$. Thus $VU' = O$, and having assumed $VU = I$, every column of $U$ has a nonzero inner product with a column of $V^*$, which implies $\col(U)\cap\col(U') = \{\bzero\}$. Since $m = \rank(VU) \leq \min(\rank(V),\rank(U))$ and $m < n$, the dimensions of $U$ and $V$ imply $\rank(U) = m$ and $\rank(V^*) = \rank(V) = m$. Thus, $U'$ has $n-m$ columns and $\left[\begin{array}{c|c} U & U' \end{array}\right]$ is nonsingular. Let $V'$ be the $(n-m)\times n$ matrix 
$V' = 
\left[\begin{array}{c|c} 
O & I 
\end{array}\right]
\left[\begin{array}{c|c} 
U & U' 
\end{array}\right]^{-1}$. 
Multiplying both sides of the defining equation for $V'$ on the right by $\left[\begin{array}{c|c} U & U' \end{array}\right]$ reveals that $V'U = O$ and $V'U' = I$. By construction,
\[
\left[\begin{array}{c} 
V \\ 
\hline 
V' 
\end{array}\right]
\left[\begin{array}{c|c} 
U & U' 
\end{array}\right] 
= 
\left[\begin{array}{c|c} 
VU & VU' \\ 
\hline 
V'U & V'U' 
\end{array}\right]
=
I. \qedhere
\]
\end{proof}

Recall that one of the goals of this paper is not only to determine which matrices are apportionable, but also to determine their apportionment constants. The remaining results of this section aid in this endeavor by providing lower bounds on apportionment constants in terms of the trace and determinant of a matrix.

\begin{lemma}\label{lem: GL lower bound}
 Let $A\in\C^{n\times n}$ be apportionable with apportionment constant $\kappa$. Then $\kappa \geq \frac{|\tr(A)|}n$.
\end{lemma}
\begin{proof}
Let $M$ be a nonsingular matrix such that $B = [b_{i,j}] = MAM^{-1}$ is uniform with common modulus $\kappa$. Then $n\kappa = \sum_{k = 1}^n |b_{k,k}| \geq \left|\sum_{k = 1}^n b_{k,k} \right| = |\tr(A)|$ and hence $\kappa \geq \frac{|\tr(A)|}n$.
\end{proof}

Let $c$ be a positive real number and let $A$ be an $n\times n$ nonsingular matrix such that each entry has modulus at most $c$. Then Hadamard's inequality \cite{BEHC2021} states that $|\det(A)|\leq c^n n^{n/2}$, with equality if and only if $\frac1{c}A$ is a complex Hadamard matrix.

\begin{lemma}\label{lem: Had lower bound}
Let $A\in \C^{n\times n}$ be apportionable with apportionment constant $\kappa$. Then $\kappa \geq n^{-1/2} |\det(A)|^{1/n}$. Moreover, if $MAM^{-1}$ is uniform with nonzero common modulus $\kappa = n^{-1/2} |\det(A)|^{1/n}$, then $\frac1{\kappa}MAM^{-1}$ is a complex Hadamard matrix.
\end{lemma}
\begin{proof}
The proof is immediate in the case $A$ is singular since it is always the case that $\kappa \geq 0$. So, suppose that $A$ is nonsingular. Let $M$ be a nonsingular matrix such that $B = MAM^{-1}$ is uniform with common modulus $\kappa$. Then Hadamard's inequality implies $|\det(A)| \leq \kappa^n n^{n/2}$ with equality if and only if $\frac1{\kappa}MAM^{-1}$ is a complex Hadamard matrix. Solving the inequality for $\kappa$ yields the desired result.
\end{proof}

When an apportionable matrix $A\in\C^{n\times n}$ is singular, Lemma~\ref{lem: Had lower bound} offers a trivial bound, making Lemma~\ref{lem: GL lower bound} potentially more useful. If $A$ is traceless ($\tr(A) = 0$), then the bound in Lemma~\ref{lem: GL lower bound} becomes trivial, and Lemma~\ref{lem: Had lower bound} may be more useful.

\section{Apportioning Nilpotent Matrices}\label{sec: nilpotent}

Recall that a square matrix $A$ is \emph{nilpotent} provided $A^k = O$ for some integer $k$. Thus, if $A$ is nilpotent, $0$ is the only eigenvalue of $A$. By \eref{eqn: Jordan outer products}, the Jordan form of an $n\times n$ nilpotent matrix $A$ can be written as $J(A) = \sum_{k=1}^{n-1} \alpha_k \be_k\be_{k+1}^\top$, where each $\alpha_k\in\{0,1\}$. The main result of this section, Theorem~\ref{thm: nilpotent}, uses this observation to establish that every nilpotent matrix $A$ is apportionable and $\ac(A) = (0,\infty)$.

The following example is provided as a guiding reference for the proof of Theorem~\ref{thm: nilpotent}.  Note that, for a nonsingular matrix $M$, $MAM^{-1}$ is uniform if and only if $MA\adj(M)$ is uniform.

\begin{example}\label{ex: nilpotent 5by5}
In this example we determine values $d_1,\ldots,d_5$ such that
\[
A =
\left[\begin{array}{ccccc}
0 & 1 & 0 & 0 & 0\\
0 & 0 & 1 & 0 & 0\\
0 & 0 & 0 &0 & 0\\
0 & 0 & 0 & 0 & 1\\
0 & 0 & 0 & 0 & 0
\end{array}\right]
\quad \text{is apportioned by} \quad
M =
\left[\begin{array}{ccccc}
1 & 0 & 0 & 0 & d_5\\
d_1 & 1 & 0 & 0 & 0\\
0 & d_2 & 1 & 0 & 0\\
0 & 0 & d_3 & 1 & 0\\
0 & 0 & 0 & d_4 & 1
\end{array}\right].
\]
Assume that $M$ is nonsingular, i.e., $\det(M) = 1 + d_1d_2d_3d_4d_5 \not= 0$. It is readily verified that 
\[
\adj(M) = 
\left[\begin{array}{ccccc}
1 & d_2d_3d_4d_5 & -d_3d_4d_5 & d_4d_5 & -d_5\\
-d_1 & 1 & d_1d_3d_4d_5 & -d_1d_4d_5 & d_1d_5\\
d_1d_2 & -d_2 & 1 & d_1d_2d_4d_5 & -d_1d_2d_5\\
-d_1d_2d_3 & d_2d_3 & -d_3 & 1 & d_1d_2d_3d_5\\
d_1d_2d_3d_4 & -d_2d_3d_4 & d_3d_4 & -d_4 & 1
\end{array}\right]
\]
and so
\begin{align*}
MA\adj(M) &= M\be_1\be_2^\top \adj(M) + M\be_2\be_3^\top \adj(M) + M\be_4\be_5^\top \adj(M)\\
&= 
\left[\begin{array}{ccccc}
-d_1 & 1 & d_1d_3d_4d_5 & -d_1d_4d_5 & d_1d_5\\
d_1d_2 - d_1^2 & -d_2 + d_1 & 1 + d_1^2d_3d_4d_5 & d_1d_2d_4d_5 - d_1^2d_4d_5 & -d_1d_2d_5 + d_1^2d_5\\
d_1d_2^2 & -d_2^2 & d_2 & d_1d_2^2d_4d_5 & -d_1d_2^2d_5\\
d_1d_2d_3d_4 & -d_2d_3d_4 & d_3d_4 & -d_4 & 1 \\
 d_1d_2d_3d_4^2 & -d_2d_3d_4^2 & d_3d_4^2 & -d_4^2 & d_4
\end{array}\right].
\end{align*}
To ensure the $(5,2)$-entry and the entries in the first row of $MA\adj(M)$ all have the same modulus, it is necessary that each $d_j$ has modulus $1$, i.e., $d_j = e^{\nu_j \ii}$ for some $\nu_j\in\R$.

Setting the modulus of the $(2,2)$-entry of $MA\adj(M)$ equal to $1$ we obtain $1 = |d_1 - d_2| =|1 - e^{(\nu_2 - \nu_1)\ii}|$. With the goal of proving Theorem~\ref{thm: nilpotent} in mind, it is worth noting that if the $(3,4)$-entry of $A$ had been $1$, the $(3,3)$-entry and $(4,4)$-entry would similarly require $1 = |d_j - d_{j+1}| = |1 - e^{(\nu_{j+1} - \nu_{j})\ii}|$ for $j = 2,3$. Thus, it is reasonable to choose $\nu_j = \tfrac{(j-1)\pi}{3}$ for $j = 1,\ldots, 4$. 

Setting the modulus of the $(2,3)$-entry of $MA\adj(M)$ equal to 1 gives 
\[
1 = |1 + d_1^2d_3d_4d_5| = \Big|1 + (d_1d_2^{-1})\prod_{j=1}^5 d_j \Big| = \Big|1 + \big(e^{(\nu_1 - \nu_2)\ii}\big)\prod_{j=1}^5 d_j \Big|.
\]
A choice of $d_5 = -\exp\Big(\Big(\tfrac{2\pi}{3} - \sum_{j=1}^4\nu_j\Big) \ii \Big) = -e^{2\pi\ii/3}$ does the job and ensures $\det(M)\not=0$. With these values of $d_1,\ldots,d_5$ we obtain
\[
MAM^{-1} = 
\frac1{1 - e^{2\pi\ii/3}}\left[\begin{array}{ccccc}
-1 & 1 & -e^{\pi\ii/3} & -e^{2\pi\ii/3} & -e^{2\pi\ii/3}\\
e^{2\pi\ii/3} & -e^{2\pi\ii/3} & -e^{2\pi\ii/3} & -e^{\pi\ii/3} & -e^{\pi\ii/3}\\
e^{2\pi\ii/3} & -e^{2\pi\ii/3} & e^{\pi\ii/3} & -e^{\pi\ii/3} & -e^{\pi\ii/3}\\
1 & -1 & -e^{2\pi\ii/3} & 1 & 1\\
-1 & 1 & e^{2\pi\ii/3} & -1 & -1
\end{array}\right],
\]
which is uniform with common modulus $\tfrac1{\sqrt{3}}$. See \cite{code} for verification of these calculations using Python \cite{code}.
\end{example}

We now utilize the techniques in Example~\ref{ex: nilpotent 5by5} to prove the main result of this section.

\begin{theorem}\label{thm: nilpotent}
Every complex nilpotent matrix $A$ is apportionable, and if $A\not=O$, then $\ac(A) = (0,\infty)$.
\end{theorem}
\begin{proof}
Throughout the proof, all index arithmetic is done modulo $n$, with indices in the residue system  $\{1,\ldots,n\}$. All empty products are equal to $1$.

Let $A\in\C^{n\times n}$ be nilpotent. Without loss of generality, assume that $A$ is in Jordan form. If $A = O$, then we are done. So, suppose $A\not= O$. Note that this implies $n\geq 2$.  By Lemma \ref{lem:pad-by-0}, we may also assume that every Jordan block of $A$ has order at least $2$. Thus, $A = \sum_{j=1}^{n} \alpha_j \be_j \be_{j+1}^\top$, where each $\alpha_k\in \{0,1\}$ with $\alpha_1 = 1$ and $\alpha_n = 0$. Note that there do not exist consecutive zeros along the superdiagonal of $A$, i.e., $\max\{\alpha_j,\alpha_{j+1}\} = 1$ for all $j = 1,\ldots,n$.

Let $P$ be the $n\times n$ permutation matrix for which $P \be_j = \be_{j+1}$ for each $\be_j$, and let $D = \diag(d_1,\ldots,d_n)$, where 
\[
d_j = \exp\left(\tfrac{(j-1)\pi}{3} \ii \right) \text{ for } j = 1,\ldots, n-1, \text{ and } d_n = \exp\Bigg(\Big(\tfrac{2\pi}{3} - \pi n - \sum_{j=1}^{n-1}\tfrac{(j-1)\pi}{3}\Big) \ii \Bigg).
\]
We show that the matrix $M = I + PD$ apportions $A$ by demonstrating $B = [b_{r,s}] = MA\adj(M)$ is uniform. Note that $M$ is indeed nonsingular since
\begin{equation}\label{eqn: det(M)}
\det(M) = 1 + (-1)^{n-1} \prod_{j=1}^n d_j = 1-\exp\left(\tfrac{2\pi}{3} \ii \right).
\end{equation}

With the goal of verifying $\adj(M)= \sum_{j=0}^{n-1} (-PD)^j$, observe that $(PD)^j$ can be written as the telescoping product $ P^j\prod_{\ell=1}^j P^{-j+\ell}D P^{j-\ell}$. Recall that $P^i\be_k = \be_{k+i}$ and notice $\be_k^\top P^i = (P^{-i} \be_k)^\top = \be_{k-i}^\top$. Since $D = \sum_{k=1}^n d_k\be_k\be_k^\top$, we have $(PD)^j =  P^j\prod_{\ell=1}^{j} \left(\sum_{k=1}^n d_k \be_{k-j+\ell}\be_{k-j+\ell}^\top\right)$. Reindexing the sum yields
\begin{equation}\label{eqn: PD to the j}
(PD)^j = P^j\prod_{\ell=1}^j \left(\sum_{k=1}^n d_{k+j-\ell} \be_{k}\be_{k}^\top\right) =  P^j\sum_{k=1}^n \left(\prod_{\ell=0}^{j-1} d_{k+\ell}\right) \be_{k}\be_{k}^\top.
\end{equation}
Observe $M\sum_{j=0}^{n-1}(-PD)^j = \sum_{j=0}^{n-1} (I + PD)(-PD)^j$ is a telescoping sum that equals $I + (-1)^{n-1}(PD)^n$. It follows from \eref{eqn: PD to the j} and \eref{eqn: det(M)} that
\[
M\sum_{j=0}^{n-1}(-PD)^j = I + (-1)^{n-1}\left(\prod_{j=1}^n d_j\right)I = \det(M)I,
\] 
proving the claim $\adj(M)= \sum_{j=0}^{n-1} (-PD)^j$. 


Further, \eref{eqn: PD to the j} implies
\[
\adj(M)
= \sum_{j=0}^{n-1} (-1)^{j} P^j \sum_{k=1}^n \left(\prod_{\ell=0}^{j-1}d_{k+\ell}\right) \be_k \be_k^\top 
= \sum_{j=0}^{n-1} (-1)^{j} \sum_{k=1}^n \left(\prod_{\ell=0}^{j-1}d_{k+\ell}\right) \be_{k+j} \be_k^\top
\]
and so the $(r,s)$-entry of $\adj(M)$ is 
\[
\be_r^\top\adj(M)\be_s = \sum_{j=0}^{n-1} (-1)^{j} \sum_{k=1}^n \left(\prod_{\ell=0}^{j-1}d_{k+\ell}\right) \be_r^\top\be_{k+j} \be_k^\top \be_s= (-1)^{t_{r,s}} \prod_{\ell=0}^{t_{r,s}-1} d_{s+\ell},
\]
where $t_{r,s} = (r-s)\mod n$; unlike indices $t_{r,s}$ is computed using the least residue system modulo $n$, i.e., $t_{r,s}\in\{0,\ldots,n-1\}$. Thus, the $(r,s)$-entry of $B$ is
\begin{align}
b_{r,s} &= \be_r^\top MA\adj(M)\be_s =  \be_r^\top (I + PD) \left(\sum_{j=1}^{n} \alpha_j \be_j\be_{j+1}^\top \right) \adj(M) \be_s \nonumber\\
&= (\be_r^\top + d_{r-1}\be_{r-1}^\top)\left( \sum_{j=1}^{n} \alpha_j \be_j \left(\be_{j+1}^\top \adj(M)\be_s\right) \right) \nonumber\\
&= \sum_{j=1}^{n} \alpha_j(\be_r^\top \be_j + d_{r-1}\be_{r-1}^\top\be_{j})(-1)^{t_{j+1,s}}\prod_{\ell=0}^{t_{j+1,s}-1} d_{s+\ell} \nonumber\\
&= \alpha_r (-1)^{t_{r+1,s}}\prod_{\ell=0}^{t_{r+1,s}-1} d_{s+\ell} + \alpha_{r-1} d_{r-1} (-1)^{t_{r,s}}\prod_{\ell=0}^{t_{r,s}-1} d_{s+\ell} \label{eqn1: nilpotent-no-nontrivial-blocks}.
\end{align}

Recall that $\min\{\alpha_{r-1},\alpha_r\} = 1$ and $|d_j| = 1$ for every $j \in \{1,\ldots, n\}$. Thus, \eref{eqn1: nilpotent-no-nontrivial-blocks} implies $|b_{r,s}| = 1$ whenever $\alpha_r = 0$ or $\alpha_{r-1} = 0$. So, suppose that $\alpha_r = \alpha_{r-1} = 1$. Note that $1 < r < n$ since $\alpha_n = 0$. We consider the cases $t_{r,s} = n-1$ and $0 \leq t_{r,s} \leq n-2$ separately. 

First, assume $t_{r,s} = n-1$. Then $t_{r+1,s} = 0$ and $s = r - (n-1) \mod n$. Since $r\notin\{1, n\}$, 
\[
d_{r-1}\prod_{\ell\not=r} d_\ell = d_{r-1}d_r^{-1}\exp\left(\left(\tfrac{2\pi}3 - \pi n \right)\ii \right) = \exp\left(\left(\tfrac{\pi}3 - \pi n\right)\ii \right).
\]
Thus, by \eref{eqn1: nilpotent-no-nontrivial-blocks}
\[
|b_{r,s}| 
= \Bigg| 1 + (-1)^{t_{r,s}}d_{r-1}\prod_{\ell=0}^{n-2} d_{s+\ell}\Bigg| 
= \Bigg|1 - \exp(n\pi\ii) d_{r-1}\prod_{\ell\not=r} d_\ell \Bigg| 
= \left|1 - \exp\left(\tfrac{\pi}3 \ii\right)\right| 
= 1.
\]

Now assume $0 < t_{r,s} < n -1$. Then $t_{r+1,s} = t_{r,s} + 1$. Since $r\notin\{1, n\}$ we have $d_r - d_{r-1} = d_r\big(1 - \exp\big(-\tfrac{\pi}{3}\ii\big)\big)$. Thus, by \eref{eqn1: nilpotent-no-nontrivial-blocks}
\[
|b_{r,s}| 
= \Bigg|  \pm \Bigg(\prod_{\ell=0}^{t_{r,s}} d_{s+\ell} - d_{r-1}\prod_{\ell=0}^{t_{r,s}-1} d_{s+\ell} \Bigg)\Bigg| 
= \Bigg| \big(d_r - d_{r-1} \big)\prod_{\ell = 0}^{t_{r,s}-1} d_\ell \Bigg| 
= \left|1 - \exp\left(-\tfrac{\pi}3 \ii\right)\right| 
= 1.
\]
Therefore, $B$ is uniform.

It remains to show that $\ac(A) = (0,\infty)$. We have already established that $A$ is apportionable with some nonzero apportionment constant $\kappa$. Thus, there exists a nonsingular matrix $M$ such that $MAM^{-1}$ is uniform with common modulus $\kappa$. Let $\kappa'\in(0,\infty)$. By Observation \ref{obs: Jordan form} $A$ and $\kappa'\kappa^{-1} A$ have the same Jordan form, i.e., there exists a nonsingular matrix $N$ such that $NAN^{-1} = \kappa'\kappa^{-1} A$. Thus, $(MN)A(MN)^{-1} = \kappa'\kappa^{-1} MAM^{-1}$ is uniform with common modulus $\kappa'$. 
\end{proof}

Having established that all nilpotent matrices are apportionable, it is reasonable to ask if Lemma~\ref{lem:pad-by-0} can be improved:  

\begin{question}\label{quest:pad-by-nil-jordan-block}
Is $A \oplus J_m(0)$ apportionable for every $m\geq1$ and every apportionable matrix $A$?
\end{question}

Thus far, the authors have only affirmed Question~\ref{quest:pad-by-nil-jordan-block} for matrices of order $n\leq 3$ (see Example~\ref{ex: 3x3 I oplus nilpotent}).

\section{Apportioning Matrices with at Most Half Rank}\label{sec: half rank}

This section's primary result, Theorem~\ref{thm: half rank apportionable}, proves that every $A\in\C^{n\times n}$ with $\rank(A)\leq \frac n2$ is apportionable. To pave the way, we first prove Theorem~\ref{thm: A oplus 0}. This theorem, while essentially equivalent to our main result, allows us to address the question: For any matrix $A$, does there exist an $m$ such that $A\oplus O_m$ is apportionable? This specific inquiry motivates Problem~\ref{prob: add zeros}. Since the core techniques in Theorem~\ref{thm: A oplus 0} can get buried in the technical details, we include Example~\ref{ex:I oplus O} as a warm-up. This example also demonstrates that the bound in Lemma ~\ref{lem: GL lower bound} is sharp.

\begin{example}\label{ex:I oplus O}
We show that $I_n\oplus O_n$ is apportionable with $\ac(I_n\oplus O_n) = \left[\tfrac12, \infty\right)$. Let $\kappa\in \big[\tfrac12,\infty \big)$ and $\zeta = \frac{1}{2} + \sqrt{\kappa^2 - \frac{1}{4}}\ii$. Observe that $\Re(\zeta)=\frac{1}{2}$ and $|\zeta|=\kappa$. For each $k \in \{1,\ldots,n\}$, define $\bu_k,\bv_k\in\C^{2n}$ via
\[
\bu_k = \zeta \sum_{j=2k}^{2n} \be_j - \overline{\zeta}\sum_{j=1}^{2k-1} \be_j = \begin{bmatrix}
    -\overline{\zeta} \\ \vdots \\ -\overline{\zeta} \\ \zeta \\ \vdots \\ \zeta
\end{bmatrix} \quad \text{ and } \quad \bv_k = \be_{2k}-\be_{2k-1} = \begin{bmatrix}
    0 \\ \vdots \\ -1 \\ 1 \\ \vdots \\ 0
\end{bmatrix}.
\]
By construction, $\bv_k^\top\bu_k = \zeta + \overline{\zeta} = 1$ for all $k \in \{1,\ldots, n\}$, and $\bv_\ell^\top\bu_k = 0$ provided $k\neq\ell$. Thus, $U = \left[\begin{array}{ccc} \bu_1 & \cdots & \bu_n \end{array}\right]$ and $V = \left[\begin{array}{ccc} \bv_1 & \cdots & \bv_n \end{array}\right]^\top$ satisfy $VU = I_n$. By Lemma \ref{lem:block-matrix-inverse}, there exists a matrix $U'\in \C^{2n\times n}$ such that $M = \left[\begin{array}{c|c} U & U' \end{array}\right]$ is nonsingular and, for $k = 1,\ldots,n$, row $k$ of $M^{-1}$ is $\bv_k^{\top}$. Observe that
\begin{align*}
    B &= M(I_n \oplus O_n)M^{-1} = M\left(\sum_{k=1}^n \be_k\be_k^\top \right)M^{-1} = \sum_{k=1}^n \left( M\be_k\right)\left(\be_k^\top M^{-1}\right) = \sum_{k=1}^n \bu_{k}\bv_{k}^\top \\
    &= \sum_{k=1}^n  \left[\begin{array}{c|c|c|c}
        O_{(2k-2)\times n} & -\bu_k & \bu_k & O_{(2n-2k)\times n}
    \end{array}\right].
\end{align*}
Since each outer product contributes two distinct uniform columns to the sum, the matrix $B$ is uniform with common modulus $|\zeta| = \kappa$. 
It follows that $\left[\frac12,\infty\right)\subseteq \ac(I_n\oplus O_n)$. Since $\frac{\tr(I_n\oplus O_n)}{2n} = \frac12$, Lemma~\ref{lem: GL lower bound} guarantees $\ac(I_n\oplus O_n) = \left[\frac12,\infty\right)$.
\end{example}

In Example~\ref{ex:I oplus O} we defined the vectors $\bu_1,\ldots,\bu_n$ and $\bv_1,\ldots,\bv_n$ so that their sum of outer products, $\sum_{k=1}^n\bu_k \bv_k^\top$, was uniform and hence $I\oplus O$ was apportionable. Verifying that this matrix was uniform was aided by the fact that the outer products in the sum had disjoint support and their nonzero entries all had the same modulus. For Theorem~\ref{thm: A oplus 0} a similar choice of vectors will be used to apportion a matrix of the form $A\oplus O$ (in Jordan form). However, to account for extra cases, their definition is more complicated, and the analogous sum is not nearly as well behaved. 

Note that, unlike Example~\ref{ex:I oplus O}, Theorem~\ref{thm: A oplus 0} does not fully determine the set of apportionment constants for the more general setting. 

\begin{theorem}\label{thm: A oplus 0}
Let $A\in \C^{n\times n}$ and assume $\rank(A)=r \geq \frac{n}2$. Then $A \oplus O_{2r-n} \in \C^{2r \times 2r}$ is apportionable and $\big(\frac{\rho(A)}{2},\infty\big) \subseteq \ac(A\oplus O_{2r-n})$. In particular, $A$ is apportionable when $\rank(A) = \frac{n}{2}$.
\end{theorem}
\begin{proof}
By Theorem~\ref{thm: nilpotent} we may assume $A$ is not nilpotent. So, $\rho(A) > 0$. Let $\kappa \in \big(\frac{\rho(A)}{2},\infty\big)$ and $B = \kappa^{-1}A$. It suffices to show that the $2r\times 2r$ matrix $B \oplus O_{2r-n}$ is apportionable and has apportionment constant $1$. Note that the eigenvalues of $B$ all have modulus strictly less than $2$. Without loss of generality, assume $B$ is in Jordan form (see \eref{eqn: Jordan outer products}):
\[
    B = \sum_{k=1}^{n} \lambda_k \be_k\be_k^\top + \sum_{k=2}^{n} \alpha_{k-1} \be_{k-1}\be_{k}^\top,
\] 
where each $\alpha_k \in \{0,1\}$ and $|\lambda_k| \geq |\lambda_{k+1}|$ for $k = 1,\ldots,n-1$.

We now use Lemma~\ref{lem:block-matrix-inverse} to construct a matrix that will be used to apportion $B \oplus O_{2r-n}$. For each $k\in\{1,\ldots,r\}$ such that $\lambda_k \neq 0$ let $\zeta_k = \frac{|\lambda_k|}{2} + \frac{\sqrt{4-|\lambda_k|^2}}{2}\ii$ and $\gamma_k = \Big(\left(\overline{\zeta_k}^2-1\right)\lambda_k\Big)^k$. Note that $|\zeta_k| = 1$ and $\Im(\zeta_k) = \frac{\sqrt{4-|\lambda_k|^2}}{2}$ which implies $\gamma_k \not= 0$ since $|\lambda_k| < 2$. For each $k\in\{1,\ldots,r\}$ such that $\lambda_k \neq 0$ define
\[
\bu_k =
\frac{1}{\gamma_{k}|\lambda_k|} 
\left(\zeta_k\sum_{j=2k}^{2r} \be_j - \overline{\zeta_k}\sum_{j=1}^{2k-1} \be_j\right) \in \C^{2r}
\qquad \text{and} \qquad 
\bv_k = 
\gamma_{k}(\be_{2k}-\be_{2k-1}) \in \C^{2r},
\]
and for each $k\in\{1,\ldots,r\}$ such that $\lambda_k = 0$ define
\[
\bu_k = e^{\ii\pi/3}\sum_{j=2k}^{2r} \be_j -e^{5\ii\pi/3}\sum_{j=1}^{2k-1} \be_j \in \C^{2r}
\qquad \text{and} \qquad
\bv_k = \be_{2k}-\be_{2k-1}\in \C^{2r}.
\]
For each $k\in \{1,\ldots, r\}$ define
\[
\widehat{\bu}_k = \sum_{j=2k+1}^{2r} \be_j - \sum_{j=1}^{2k} \be_j \in \C^{2r}.
\]
For $k \in \{r + 1,\ldots, 2r\}$ let $\bu_k = \widehat{\bu}_{k - r}$. It is readily verified that $\bv_k^\top\bu_k = 1$ and $\bv_k^\top\bu_\ell =0$ for $\ell \neq k$, i.e.,
\[
\left[\begin{array}{c} \bv_1^\top\\\hline \vdots \\\hline\\[-1em] \bv_r^\top \end{array}\right]
\left[\begin{array}{c|c|c} \bu_1 & \cdots & \bu_r \end{array}\right]
 = I_r,
\]
Moreover, $\{\widehat{\bu}_1,\ldots,\widehat{\bu}_r\}$ is a basis for $\spn\{\bv_1,\ldots,\bv_r\}^\perp$. Thus, by Lemma \ref{lem:block-matrix-inverse}, the $2r\times 2r$ matrix $M = \left[\begin{array}{c|c|c} \bu_1 & \cdots & \bu_{2r} \end{array}\right]$ is nonsingular and the $k$-th row of $M^{-1}$ is $\bv_k^\top$ for $k\in\{1,\ldots,r\}$.

To obtain a matrix that apportions $B \oplus O_{2r-n}$, we must permute the columns of $M$. Let 
\[
\Omega = \{1\}\cup\{\ell \in \{2,\ldots,n\}: \lambda_\ell \neq 0 \text { or } \alpha_{\ell-1} \neq 0\},
\]
i.e., $\Omega$ is the set of indices corresponding to nonzero columns of B. Note that $|\Omega| = r$ since $B$ is not nilpotent. Let $\varphi$ be any bijection on $\{1,\ldots, 2r\}$ such that $\varphi|_\Omega$ maps the $j$-th smallest element of $\Omega$ to $j$. So, $\{\varphi(k): k\in\Omega\} = \{1,\ldots, r\}$ and $\varphi|_{\{k: \lambda_k \neq 0\}}$ is the identity map. Let $P$ be the $2r\times 2r$ permutation matrix such that $P\be_j = \be_{\varphi(j)}$ for $j\in\{1,\ldots,2r\}$. We claim that $MP$ apportions $B\oplus O_{2r-n}$ with common modulus $1$.

Let $C = (MP)(B \oplus O_{2r-n})(MP)^{-1}$. Direct computation yields
\begin{align*}
C &= MP\left(\sum_{k=1}^{n} \lambda_k \be_k\be_k^\top + \sum_{k=2}^{n} \alpha_{k-1} \be_{k-1}\be_{k}^\top \right)P^\top M^{-1}\\
&= M\left(\sum_{k\in\Omega} \lambda_k \be_{\varphi(k)}\be_{\varphi(k)}^\top + \sum_{k\in\Omega\setminus\{1\}} \alpha_{k-1} \be_{\varphi(k-1)}\be_{\varphi(k)}^\top \right) M^{-1}\\
&= \sum_{k\in\Omega} \lambda_k \bu_{\varphi(k)}\bv_{\varphi(k)}^\top + \sum_{k\in\Omega\setminus\{1\}} \alpha_{k-1} \bu_{\varphi(k-1)}\bv_{\varphi(k)}^\top.
\end{align*}
Let $i,j \in \{1,\ldots,2r\}$. Since every pair of vectors taken from $\{\bv_1,\ldots,\bv_r\}$ have disjoint support and $\bigcup_{k=1}^r \supp(\bv_k) = \{1,\ldots,2r\}$, there is a unique $\ell\in \{1,\ldots,r\}$ such that $\bv_{\ell}^\top\be_j\not=0$ (specifically, $\ell = \big\lceil \frac j2 \big\rceil)$.  Note that if $\ell = \varphi(k)$, then $k-1 = \varphi^{-1}(\ell) - 1$. For brevity, write $\ell' = \varphi^{-1}(\ell) - 1$ and let $c_{i,j}$ denote the $(i,j)$-entry of $C$. Then
$c_{i,j} = \lambda_1 \be_i^\top \bu_1 \bv_1^\top\be_j$ for $j\in\{1, 2\}$. For $j\geq 3$, since $\lambda_\ell = \lambda_{\varphi^{-1}(\ell)}$,
\[
c_{i,j} = \be_i^\top C \be_j
= \lambda_{\ell} \be_i^\top \bu_{\ell}\bv_{\ell}^\top\be_j + \alpha_{\ell'} \be_i^\top\bu_{\varphi(\ell')}\bv_{\ell}^\top\be_j.
\]
There are three cases to consider: $\lambda_\ell = 0$ and $\alpha_{\ell'} = 1$, $\lambda_\ell \neq 0$ and $\alpha_{\ell'} = 0$, and $\lambda_\ell \neq 0$ and $\alpha_{\ell'} = 1$.
    
Assume that $\lambda_\ell = 0$ and $\alpha_{\ell'} = 1$. Then the $j$-th entry of $\bv_\ell$ is equal to $\pm 1$. Note that $\bu_{\varphi(\ell')} = \bu_k$ or $\bu_{\varphi(\ell')} = \widehat{\bu}_k$ for some $k\in\{1,\ldots,r\}$. So, the $i$-th entry of $\bu_{\varphi(\ell')}$ is in the set $\{e^{i\pi/3},-e^{5i\pi/3}, 1, -1\}$. Thus $|c_{i,j}| = \left|\alpha_{\ell'} \be_i^\top\bu_{\ell'}\bv_{\ell}^\top\be_j \right| = 1$. 

Assume $\lambda_\ell \neq 0$ and $\alpha_{\ell'} = 0$. Then the $i$-th entry of $\bu_\ell$ is in $\left\{\frac{\zeta_\ell}{\gamma_{\ell}|\lambda_\ell|},\frac{-\overline{\zeta_\ell}}{\gamma_{\ell}|\lambda_\ell|}\right\}$ and the $j$-th entry of $\bv_\ell$ is equal to $\pm \gamma_\ell$. Thus $|c_{i,j}| = \left|\lambda_{\ell} \be_i^\top\bu_{\ell}\bv_{\ell}^\top\be_j \right| = |\zeta_\ell| = 1$.

Assume $\lambda_\ell \neq 0$ and $\alpha_{\ell'} = 1$. Then the $j$-th entry of $\bv_\ell$ is equal to $\pm\gamma_\ell$. Observe that $\varphi(\ell) = \ell$ since $\lambda_\ell \neq 0$ and $|\lambda_k| \geq |\lambda_{k+1}|$ for all $k\in \{1,\ldots,n-1\}$. Thus, $\ell' = \ell - 1$ and $\varphi(\ell') = \ell - 1$. Since $\alpha_{\ell-1}=1$, we have that $\lambda_{\ell-1}$ and $\lambda_\ell$ are in the same Jordan block. So, $\lambda_{\ell-1}=\lambda_\ell$ which implies that $\zeta_{\ell-1} = \zeta_{\ell}$. Let $x_i$ denote the $i$-th entry of $\bu_\ell$ and $y_i$ denote the $i$-th entry of $\bu_{\ell-1}$. Then
\[
|c_{i,j}| = |\lambda_{\ell} \be_i^\top \bu_{\ell}\bv_{\ell}^\top\be_j + \be_i^\top\bu_{\ell - 1}\bv_{\ell}^\top\be_j| = |\lambda_{\ell} x_i(\pm\gamma_\ell) + y_i(\pm\gamma_\ell)|,
\]
and
\[
x_i = 
\begin{cases}
\frac{-\overline{\zeta_{\ell}}}{\gamma_{\ell}|\lambda_{\ell}|} & \text{if $i \leq 2\ell-1$}, \\
\frac{\zeta_{\ell}}{\gamma_{\ell}|\lambda_{\ell}|} & \text{if $i \geq 2\ell$}, 
\end{cases}
\quad \text{ and } \quad 
y_i = 
\begin{cases}
\frac{-\overline{\zeta_{\ell}}}{\gamma_{\ell-1}|\lambda_{\ell}|} & \text{if $i \leq 2\ell-3$}, \\
\frac{\zeta_{\ell}}{\gamma_{\ell-1}|\lambda_{\ell}|} & \text{if $i \geq 2\ell-2$}.
\end{cases}
\]
Suppose $i \leq 2\ell-3$ or $i \geq 2\ell$. Then either $x_i = \frac{-\overline{\zeta_{\ell}}}{\gamma_{\ell}|\lambda_{\ell}|}$ and $y_i = \frac{-\overline{\zeta_{\ell}}}{\gamma_{\ell-1}|\lambda_{\ell}|}$, or $x_i = \frac{\zeta_{\ell}}{\gamma_{\ell}|\lambda_{\ell}|}$ and  $y_i = \frac{\zeta_{\ell}}{\gamma_{\ell-1}|\lambda_{\ell}|}$. Observe that
\[
\frac{\gamma_{\ell}}{\gamma_{\ell-1}} = \frac{\Big[\left(\overline{\zeta_{\ell}}^2-1\right)\lambda_\ell\Big]^{\ell}}{\Big[\left(\overline{\zeta_\ell}^2-1\right)\lambda_\ell\Big]^{\ell-1}} = \left(\overline{\zeta_\ell}^2-1\right)\lambda_\ell.
\] 
Thus,
\[
|c_{i,j}| = |\lambda_{\ell} x_i(\pm\gamma_\ell) + y_i(\pm\gamma_\ell)|
= |\zeta_\ell|\left|\frac{\lambda_\ell}{|\lambda_\ell|}\right|\left|1 + \frac{\gamma_\ell}{\gamma_{\ell-1}\lambda_\ell} \right| = \left|1 + \left(\overline{\zeta_\ell}^2-1\right)\right| = 1.
\]
Suppose $2\ell-2 \leq i \leq 2\ell-1$.  Then $x_i = \frac{-\overline{\zeta_{\ell}}}{\gamma_{\ell}|\lambda_{\ell}|}$ and $y_i = \frac{\zeta_{\ell}}{\gamma_{\ell-1}|\lambda_{\ell}|}$. So,
\[
|c_{i,j}| = |\lambda_{\ell} x_i(\pm\gamma_\ell) + y_i(\pm\gamma_\ell)|
= |\zeta_\ell|\left|\frac{\lambda_\ell}{|\lambda_\ell|}\right|\left|\frac{-\overline{\zeta_{\ell}}}{\zeta_\ell} + \frac{\gamma_\ell}{\gamma_{\ell-1}\lambda_\ell} \right| 
= \left|-\overline{\zeta_\ell}^2 + \left(\overline{\zeta_\ell}^2-1\right)\right| 
= 1.
\]
Therefore, $C$ is uniform with common modulus $1$. 
\end{proof}

Theorem~\ref{thm: half rank apportionable} follows from Theorem~\ref{thm: A oplus 0} by utilizing some basic facts about Jordan form. Recall that if the rank of $A\in\C^{n\times n}$ is $r$, then the zero eigenvalue has geometric multiplicity $n - r$ which means $J(A)$ contains precisely $n-r$ Jordan blocks corresponding to $0$. Thus, $\rank(A) < \frac n2$ implies the number of Jordan blocks of $J(A)$ corresponding to $0$ is strictly greater than $\frac n2$ and so $J(A)$ contains at least one Jordan block of size $1$ corresponding to $0$.

\begin{observation}\label{obs:rank < half -> A = B oplus 0}
Suppose $A \in \C^{n\times n}$ with $\rank(A)<\frac{n}2$. Then there exists a $B \in \C^{n-1\times n-1}$ such that $A$ is similar to $B \oplus [0]$. 
\end{observation}

\begin{theorem}\label{thm: half rank apportionable}
Every $A \in \C^{n\times n}$ with $\operatorname{rank}(A) \leq \frac{n}2$ is apportionable and $\big(\frac{\rho(A)}{2},\infty\big) \subseteq \ac(A)$.
\end{theorem}

\begin{proof}
Since $\rank(A) \leq \frac n2$, there exists an $m\in \{0,\ldots,n\}$ such that $2(\rank(A)) = n - m$. Repeated applications of Observation~\ref{obs:rank < half -> A = B oplus 0} imply that $A$ is similar to a matrix of the form $B\oplus O_m$ with $\rank(B) = \frac{n-m}2$, i.e., its rank is half its order. Since $\rho(B) = \rho(A)$, Theorem \ref{thm: A oplus 0} implies $B$ is apportionable with $\big(\frac{\rho(A)}{2},\infty\big) \subseteq \ac(B)$. Therefore, $A$ is apportionable with $\big(\frac{\rho(A)}{2},\infty\big) \subseteq \ac(A)$ by Lemma \ref{lem:pad-by-0}.
\end{proof}

Define $\varsigma$ on the set of all square matrices via $\varsigma(A) = \min\{m : A\oplus O_m \text{ is apportionable}\}$. Note Theorem~\ref{thm: half rank apportionable} guarantees $\varsigma(A) = 0$ for all $A\in\C^{n\times n}$ with $\rank(A)\leq \frac n2$, and Theorem~\ref{thm: A oplus 0} guarantees $\varsigma(A)\leq 2\rank(A) - n \leq n$ for all matrices $A\in\C^{n\times n}$ with $\rank(A)\geq \frac n2$. Examples of nonsingular matrices $A$ in which $\varsigma(A) = n$ and $0 < \varsigma(A) < n$ are discussed after Theorem \ref{thm: I oplus lambda}. This naturally leads to the following problem.

\begin{problem}\label{prob: add zeros}
Determine $\varsigma(A)$ for all $A\in\C^{n\times n}$.
\end{problem}

\section{Small Order Apportionable Matrices}\label{sec: small order}

This section addresses apportionment for small order matrices. We first present a few tools in Section~\ref{sub:Small Tools}. These include an apportionment constant result for rank one matrices (Proposition \ref{prop: rank 1 apport}) and an apportionment result on perturbations of the identity (Theorem \ref{thm: I oplus lambda}). Section ~\ref{sec: order 2} summarizes all results for $2\times 2$ matrices, the characterization of which was initiated in \cite{CCGH2024} and is completed by Theorem~\ref{thm: 2x2 distinct evals}. As a proof of concept, a Python function in \cite{code} determines if a $2\times 2$ matrix is apportionable and provides an apportioning matrix if one exists. Section~\ref{sec: order 3} initiates the classification of all $3\times 3$ matrices.

\subsection{Useful Tools}\label{sub:Small Tools}

We note that \cite[Example 8.1]{CCGH2024} indirectly establishes $\ac(A) = \left[\frac{|\lambda|}{2},\infty\right)$ for $A = \diag(\lambda,0)$. The following proposition generalizes this result by establishing the set of apportionment constants for rank one matrices of any order. We first require the following lemma.

\begin{lemma}\label{lem: spiral sum} 
For any integer $n\geq 2$ and real number $r\geq \frac{1}{n}$ there exist real numbers $\theta_1,\ldots,\theta_n$ such that $r\sum_{j=1}^n e^{\ii\theta_j} = 1$.
\end{lemma}
\begin{proof}
Fix an integer $n\geq 2$ and let $r\geq \frac 1n$ be a real number. Consider the continuous real function $f$ 
given by $f(\theta) = \left|\sum_{j=1}^{n} e^{\ii j \theta }\right|$. Since $f(0) = n$ and $f\parens{\frac{2\pi}{n}} = 0$, there exists a real value $\rho\in \left[0,\frac{2\pi}{n}\right]$ such that $f(\rho) = \frac1r$. Let $\alpha$ be any real number satisfying $e^{\ii\alpha} f(\rho) = \sum_{j=1}^{n} e^{\ii j\rho}$. For each $j\in \{1,\ldots,n\}$ take $\theta_j = j\rho - \alpha$ so that $r\sum_{j=1}^n e^{\ii\theta_j} = r e^{-\ii\alpha} \sum_{j=1}^{n} e^{\ii j\rho}  = r f(\rho) = 1$, as desired.
\end{proof}

The following proposition provides another example for which Lemma~\ref{lem: GL lower bound} is sharp.

\begin{proposition}\label{prop: rank 1 apport} 
Let $A\in\C^{n\times n}$ and assume $\rank(A) = 1$. Then $A$ is apportionable and
\[
\ac(A) =
\begin{cases}
(0,\infty) & \text{if $A$ is nilpotent}\\
\left[\frac{\rho(A)}n, \infty \right) & \text{otherwise.}
\end{cases}
\]
\end{proposition}
\begin{proof}
If $A$ is nilpotent, then we are done by Theorem~\ref{thm: nilpotent}. So, suppose that $A$ is not nilpotent. Up to similarity $A = \diag(\lambda, 0,\ldots,0)$, where $|\lambda| = \rho(A)\not=0$. Let $\kappa\in \left[\frac{\rho(A)}n, \infty \right)$ and $r = \frac{\kappa}{\rho(A)}$. Observe that $r\geq \frac1n$. By Lemma~\ref{lem: spiral sum} there exist real numbers $\theta_1,\ldots,\theta_n$ such that $r\sum_{j=1}^n e^{\ii\theta_j} = 1$.
Let $\bu = \left[\begin{array}{ccc}
1 & \cdots & 
1 \end{array}\right]^\top$ and $\bv = \left[\begin{array}{ccc}
r e^{\ii \theta_1} & \cdots & 
r e^{\ii \theta_n} \end{array}\right]^\top$ so that $\bv^\top\bu = 1$. By Lemma~\ref{lem:block-matrix-inverse} there exists a nonsingular matrix $M$ such that the first column of $M$ is $\bu$ and the first row of $M^{-1}$ is $\bv$. Then the $(i,j)$-entry of
$MAM^{-1} = \lambda\bu\bv^\top$ is $\lambda r e^{\ii \theta_j}$ which has modulus $\kappa$. Thus, $MAM^{-1}$ is uniform and $\left[\frac{\rho(A)}{n},\infty\right) \subseteq \ac(A)$. Equality follows from Lemma \ref{lem: GL lower bound}.
\end{proof}

A \emph{rank one perturbation} of a matrix $A$ is a matrix of the form $A + B$, where $B$ is a rank one matrix. Since every rank one matrix $B$ of order $n\geq 2$ is similar to $E_{1,2}$ or $\lambda E_{n,n}$ for some $\lambda\in\C$, we have the following observation.

\begin{observation}\label{obs: rank one pert} 
Let $n\geq 2$ and $A\in \C^{n\times n}$ be a rank one perturbation of $I$. Then $A$ is similar to $I + E_{1,2}$ or $I\oplus [\lambda]$ for some $\lambda\in\C$.
\end{observation}

\begin{theorem}\label{thm: I oplus lambda}
Let $A\in\C^{n\times n}$ be a rank one perturbation of the identity for some $n\geq 3$. Then $A$ is apportionable if and only if $A$ is similar to $I\oplus [\lambda]$ for some $\lambda\in\C$ such that $\Re(\lambda) = 1 - \frac{n}2$. Furthermore, if $A$ is apportionable, then the set of apportionment constants of $A$ is  
\[
\ac(A) =
\begin{cases}
\left[\frac12,\infty\right) & \text{if $n$ is even and $\Im(\lambda) = 0$},\\[5pt]
\left\{\sqrt{\frac{\Im(\lambda)^2}{(n-2s)^2} +\frac14} : s = 0,\ldots, \left\lfloor \frac {n-1}2 \right\rfloor \right\} & \text{otherwise.}
\end{cases}
\]
\end{theorem}
\begin{proof}
We begin with the backwards direction. It suffices to assume $A = I\oplus [\lambda]$ for some $\lambda\in\C$ such that $\Re(\lambda) = 1 - \frac n2$. Let $F$ be the $n\times n$ \emph{discrete Fourier transform matrix}, i.e., the unitary matrix whose $(j,k)$-entry is $\frac1{\sqrt{n}}\omega^{(j-1)(k-1)}$, where $\omega =  e^{-2\pi \ii/ n}$. Then
\[
FAF^* = F(I - (1 - \lambda)\be_n\be_n^\top)F^* = I - \Big(\frac n2 - \Im(\lambda) \ii\Big)F\be_n \be_n^\top F^*.
\]
Observe that the diagonal entries of $F\be_n \be_n^\top F^*$ are equal to $\frac1n$ and the off diagonal entries have modulus $\frac1n$. Thus, the diagonal entries of $FAF^*$ are equal to $\frac 12 + \frac {\Im(\lambda)}n \ii$ and the off-diagonal entries of $FAF^*$ have modulus $\left|\frac12 - \frac {\Im(\lambda)}n \ii\right|$. Therefore, $A$ is apportionable.

Now suppose $A$ is apportionable. By Observation~\ref{obs: rank one pert} we may assume, without loss of generality, that $A = I + E_{1,2}$ or $A = I\oplus[\lambda]$ for some $\lambda\in \C$. We consider both cases simultaneously by observing $A = I + \bu\bv^\top$ for some $\bu,\bv\in \C^n$. Since $A$ is apportionable, there exists a nonsingular matrix $M$ and a positive real number $\kappa$ such that 
$B = MAM^{-1} = I + M\bu \bv^\top M^{-1}$ 
is uniform with common modulus $\kappa$. Let $x_i$ and $y_i$ denote the $i$-th entry of $M\bu$ and $\bv^\top M^{-1}$, respectively. Since $B$ is uniform, $\kappa = |x_i||y_j|$ for all $i\not=j$. Thus, for distinct $i,j,k\in\{1,\ldots,n\}$ (which exist since $n\geq 3$), we have $|x_i||y_k| = |x_j||y_k|$ and $|y_i||x_k| = |y_j||x_k|$, which implies $|x_i| = |x_j|$ and $|y_i| = |y_j|$. Thus, $M\bu$ and $\bv^\top M^{-1}$ are uniform vectors. It follows that both $B = [b_{i,j}]$ and $I - B = M\bu \bv^\top M^{-1}$ are uniform with common modulus $\kappa$. Comparing the diagonal entries of $I - B$ and $B$ we obtain $|b_{j,j}| = |1 - b_{j,j}|$ and so $\Re(b_{j,j}) = \frac12$ for $j\in \{1,\ldots,n\}$. Consequently, if $A = I + E_{1,2}$, then $n = \tr(A) = \tr(B) = \frac n2 + z\ii$ for some $z\in\R$, a contradiction. Thus, $A = I\oplus [\lambda]$. Since $\Re(b_{j,j}) = \frac12$ and $|b_{j,j}| = \kappa$ for $j\in \{1,\ldots,n\}$, each diagonal entry of $B$ is of the form $b_{j,j} = \frac12 \pm \sqrt{\kappa^2 - \frac14}\ \ii$ and
\begin{equation}\label{eqn2: I oplus lambda}
\kappa \in \left[\tfrac12,\infty \right).
\end{equation}
So,
\begin{equation}\label{eqn3: I oplus lambda}
n - 1 + \lambda = \tr(A) = \tr(B) = \frac{n}{2} + (n -2s)\sqrt{\kappa^2 - \frac14}\ \ii,
\end{equation}
where $s\in\{0,\ldots, n\}$ is the number of indices $j$ such that $\Im(b_{jj}) < 0$. Solving \eref{eqn3: I oplus lambda} for $\lambda$ yields $\Re(\lambda) = 1-\frac{n}{2}$, as desired. 

To complete the proof we must determine the apportionment constants of $A$. Assume that $A$ is apportionable. By the preceding argument, we may assume $A = I \oplus [\lambda]$ for some $\lambda\in\C$ such that $\Re(\lambda) = 1 - \frac n2$. For each $t\in\big[\frac12,\infty\big)$ and $r\in\{0,1,\ldots, n\}$ define the matrix
\[
M_r(t) = 
\left[
\begin{array}{ccccc}
-1 & \cdots & -1 & -1 & w_{1}(t)\\[7pt]
0 & \cdots & 0 & 1 & w_{2}(t)\\
\vdots & \iddots & 1 & 0 & w_{3}(t)\\
0 & \iddots & \iddots & \vdots & \vdots\\[3pt]
1 & 0 & \cdots & 0 & w_{n}(t)
\end{array}\right],
\]
where $w_{1}(t),\ldots, w_{r}(t) = \frac{\frac12 + \sqrt{t^2 - \frac14}\ \ii}{1-\lambda}$ and $w_{r+1}(t),\ldots,w_{n}(t) = \frac{\frac12 - \sqrt{t^2 - \frac14}\ \ii }{1-\lambda}$. Note that $M_r(t)$ is nonsingular if and only if the first row of $M_r(t)$ is not in the span of the remaining rows of $M_r(t)$ 
if and only if $\sum_{j=1}^{n} w_j(t) \not= 0$.

First, suppose that $n$ is even and $\Im(\lambda) = 0$. Let $s = \frac n2$ and $\kappa \in \big[\frac12,\infty\big)$ and consider the matrix $M_s(\kappa)$. Observe that $\sum_{j=1}^{n} w_j(\kappa) = 1$, and hence $M_s(\kappa)$ is nonsingular. Let $\mathds{1}$ be the all ones vector. Since $\sum_{j=1}^{n} w_j(\kappa) = 1$ we have $\mathds{1}^\top M_s(\kappa) = \be_n^\top$. Ergo, $\be_n^\top M_s^{-1}(\kappa) = \mathds{1}^\top$ and so
\[
M_s(\kappa) A M_s^{-1}(\kappa) = M_s(\kappa)(I - (1 - \lambda)\be_n\be_n^\top)M_s^{-1}(\kappa) = I - \frac n2 
\left[\begin{array}{ccc}
w_1(\kappa) & \cdots & w_1(\kappa)\\
\vdots &  & \vdots\\
w_n(\kappa) & \cdots & w_n(\kappa)
\end{array}\right].
\]
Since $1 - \lambda = \frac n2$ and $\kappa \geq \frac12$, the diagonal entries of $M_s(\kappa) A M_s^{-1}(\kappa)$ have modulus
\[
\left|1 - \frac n2 \left(\frac{\frac12 \pm \sqrt{\kappa^2 - \frac14}\ \ii }{1-\lambda}\right)\right| 
= \left|1 - \frac n2 \left(\frac{\frac12 \pm \sqrt{\kappa^2 - \frac14}\ \ii}{\frac n2}\right)\right|
= \left|\frac12 \pm \sqrt{\kappa^2 - \frac14}\ \ii\right| = \kappa.
\]
A similar calculation shows that the off-diagonal entries also have modulus $\kappa$. By \eref{eqn2: I oplus lambda}, $\ac(A) = [\frac12,\infty)$.

Now suppose that $n$ is odd or $\Im(\lambda) \not= 0$. Let $M$ be a nonsingular matrix that apportions $A$ with apportionment constant $\kappa$. By \eref{eqn3: I oplus lambda} there exists an $s\in\{0,\ldots,n\}$ such that $\Im(\lambda) = (n -2s)\sqrt{\kappa^2 - \frac14}$. Observe that $n - 2s \not= 0$ since $n$ is odd or $\Im(\lambda) \not= 0$. Solving for $\kappa$ yields $\kappa = \sqrt{\frac{\Im(\lambda)^2}{(n-2s)^2} +\frac14}$. Since 
$\left\{(n-2s)^2 : s = 0,\ldots, n\right\} = \left\{(n-2s)^2 : s = 0,\ldots, \left\lfloor\frac{n-1}2\right\rfloor \right\}$, we have
\[
\ac(A) \subseteq \left\{\sqrt{\frac{\Im(\lambda)^2}{(n-2s)^2} +\frac14} : s = 0,\ldots, \left\lfloor \frac {n-1}2 \right\rfloor \right\}.
\]
Let $s \in \big\{0,\ldots, \left\lfloor \frac {n-1}2 \right\rfloor\big\}$ and $\kappa_s = \sqrt{\frac{\Im(\lambda)^2}{(n-2s)^2} +\frac14}$. Similar to the argument above, $M_s(\kappa_s)$  is nonsingular and $\mathds{1}^\top M_s(\kappa_s) =\be_n^\top$ implying that $M_s(\kappa_s) A M_s^{-1}(\kappa_s)$ is uniform with common modulus $\kappa_s$. Thus, $\ac(A) = \left\{\sqrt{\frac{\Im(\lambda)^2}{(n-2s)^2} +\frac14} : s = 0,\ldots, \left\lfloor \frac {n-1}2 \right\rfloor\right\}$.
\end{proof}

\begin{remark}
Let $A = I_{n-1} \oplus \left[\, \lambda \,\right]$ such that $\Re(\lambda) = 1 - \frac{n}{2}$. The proof of Theorem~\ref{thm: I oplus lambda} shows that $A$ is $\U$-apportionable.
\end{remark}

Theorem~\ref{thm: I oplus lambda} provides the only known nontrivial lower bound on $\varsigma(I)$ (c.f.\ Problem~\ref{prob: add zeros}). Specifically, this result implies $\varsigma(I_n) \geq 2$ for $n \geq 2$. Combining this with Theorem \ref{thm: A oplus 0} gives $2 \leq \varsigma(I_n) \leq n$. Thus, $\varsigma(I_2) = 2$. It is unknown whether $\varsigma(I_3)$ is $2$ or $3$.

\subsection{Classification of Apportionable Matrices of order 2}\label{sec: order 2}

Before presenting the complete characterization of apportionable $2\times 2$ matrices, we summarize prior results. Clearly, the zero matrix $O$ is uniform and $\ac(O) = \{0\}$. By \cite[Theorem 2.2]{CCGH2024} all rank 1 matrices are U-apportionable and hence apportionable. If $A\in\C^{2\times 2}$ has a nonzero repeated eigenvalue, then $A$ is not apportionable by \cite[Theorem 7.4]{CCGH2024}. This leaves the case of nonsingular $2\times 2$ matrices with distinct eigenvalues. By \cite[Theorem 7.6]{CCGH2024}, if  $A\in\C^{2\times 2}$ is nonsingular and has real eigeinvalues, then $A$ is apportionable if and only if $\spec(A) = \{\lambda, -\lambda\}$. 

Table~\ref{tab: 2 by 2} presents which $2\times 2$ matrices are apportionable in terms of their Jordan form and includes new results determining their sets of apportionment constants.

\begin{table}[h!]
\begin{center}
\renewcommand{\arraystretch}{1.8}
\begin{tabular}{|c||*6{>{\renewcommand{\arraystretch}{1}}c|}}
\hline
Matrix $A$ & 
$\left[\begin{array}{cc} 0 & 0\\ 0 & 0 \end{array}\right]$
& 
$\left[\begin{array}{cc} 0 & 1\\ 0 & 0 \end{array}\right]$
& 
$\left[\begin{array}{cc} \lambda_1 & 0\\ 0 & 0 \end{array}\right]$
& 
$\left[\begin{array}{cc} \lambda_1 & 0\\ 0 & \lambda_1 \end{array}\right]$
& 
$\left[\begin{array}{cc} \lambda_1 & 1\\ 0 & \lambda_1 \end{array}\right]$
& 
$\left[\begin{array}{cc} \lambda_1 & 0\\ 0 & \lambda_2 \end{array}\right]$\\[3pt]
\hline
Apportionable & Yes & Yes & Yes & No & No & Depends on $\lambda_1, \lambda_2$
\\[3pt]
\hline
$\ac(A)$ & $\{0\}$ & $(0,\infty)$ & $\Big[\frac{|\lambda_1|}2, \infty \Big)$ & $\emptyset$ & $\emptyset$ & Depends on $\lambda_1, \lambda_2$
\\[3pt]
\hline
Proof & $A = O$ & Thm.~\ref{thm: nilpotent} & Prop.~\ref{prop: rank 1 apport} & $A = \lambda_1 I$ & \cite[Thm.\ 7.4]{CCGH2024} & Thm.~\ref{thm: 2x2 distinct evals}
\\[3pt]
\hline
\end{tabular}
\end{center}
\caption{Apportionment constants for $2\times 2$ matrices in Jordan form ($\lambda_1$ and $\lambda_2$ are distinct and nonzero).}\label{tab: 2 by 2}
\end{table}

The following theorem provides exact conditions for when $2\times 2$ matrices with rank $2$ and distinct eigenvalues are apportionable, along with all possible apportionment constants. To simplify calculations we make the following observation and leave its verification as an exercise to the reader.

\begin{observation}\label{obs: 2x2 distinct evals} 
Let $z_1, z_2, z_3 \in \C$. Then
\[
\frac{|z_1 - z_2|^2 - |z_3 - z_2|^2}2 = \frac{|z_1|^2 - |z_3|^2}{2} - \Re(z_1\overline{z_2}) + \Re(z_3\overline{z_2}).
\]
\end{observation}


\begin{theorem}\label{thm: 2x2 distinct evals}
Let $A\in \mathbb{C}^{2\times 2}$ have distinct nonzero eigenvalues $\lambda_1$ and $\lambda_2$, and let $\gamma = \frac{\lambda_2 + \lambda_1}{\lambda_2 - \lambda_1}$. Then $A$ is apportionable if and only if $\gamma = 0$ or $\Re(\gamma^2) < |\gamma|^4 \leq 1$. Moreover, if $A$ is apportionable, then
\[
\ac(A) = \begin{cases}
\left[\frac{\rho(A)}{\sqrt{2}},\infty\right) & \text{ if }\gamma = 0 \\[5pt]
\set{\left|\frac{\lambda_1 + \lambda_2}{2}\right|\sqrt{1+\frac{1-|\gamma|^4}{2(|\gamma|^4 - \Re(\gamma^2))}}} & \text{ otherwise.}
\end{cases}
\]
\end{theorem}
\begin{proof}
Without loss of generality $A = \diag(\lambda_1,\lambda_2)$. Assume that $A$ is apportionable. Then 
there exists a matrix 
$M=\begin{bmatrix}
a & b \\ c & d
\end{bmatrix}$ 
such that $\det(M) = ad - bc = 1$ and $B = MAM^{-1}$ is uniform. Let $\omega = 2bc +1$. 
Utilizing the substitution $ad = 1 + bc$ it is readily verified via direct computation that
\begin{equation}\label{eqn: mtx B}
B = (\lambda_2-\lambda_1)
\left[\begin{array}{cc}
\frac{\gamma-\omega}{2} & ab \\
-cd & \frac{\gamma + \omega}{2} 
\end{array}\right].
\end{equation}
Since $B$ is uniform $|\gamma - \omega| = |\gamma + \omega|$. Further, $|\gamma^2 - \omega^2| = |1 - \omega^2|$ since
\[
|\gamma^2 - \omega^2| = |\gamma - \omega||\gamma +\omega| = 4|ab||-cd| = 4|adbc|= 4|(1 + bc)bc| 
= |-4(bc)^2 - 4bc| = |1 - \omega^2|.
\]
For brevity, write $x = \Re(\omega)$ and $y = \Im(\omega)$. By Observation~\ref{obs: 2x2 distinct evals},
\[
0 
= \frac{|\gamma+\omega|^2 - |\gamma - \omega|^2}4
= \Re(\gamma\overline{\omega}) = \Re(\gamma)x + \Im(\gamma)y
\]
and (note $\Re(\gamma\overline{\omega}) = 0$)
%
\begin{align*}
0 &= \frac{|1 - \omega^2|^2 - |\gamma^2 - \omega^2|^2}2 \\
&= \frac{1 - |\gamma|^4}{2} - \Re(\omega^2) + \Re\big((\gamma\overline{\omega})^2\big) \\
&= \frac{1 - |\gamma|^4}{2} - x^2 + y^2 + 2\Re(\gamma\overline{\omega})^2 - |\gamma \overline{\omega}|^2 \\
&= \frac{1 - |\gamma|^4}{2} - x^2 + y^2 - |\gamma \overline{w}|^2\\
&= \frac{1 - |\gamma|^4}{2} - x^2 + y^2 - |\gamma|^2 (x^2 + y^2) \\
&= \frac{1 - |\gamma|^4}{2} - (1 + |\gamma|^2)x^2 + (1 - |\gamma|^2)y^2.
\end{align*}
Thus far we have shown $x = \Re(\omega)$ and $y = \Im(\omega)$ satisfy the system 
\begin{subequations}
\label{system: x and y}
\begin{align}
\Re(\gamma)x + \Im(\gamma)y &= 0 \label{eqn: line}\\
(1 + |\gamma|^2)x^2 - (1 - |\gamma|^2)y^2 &= \frac{1 - |\gamma|^4}{2}. \label{eqn: hyperbola}
\end{align}
\end{subequations}

We finish establishing that $\gamma = 0$ or $\Re(\gamma^2) < |\gamma|^4 \leq 1$ by considering the cases $|\gamma|>1$, $\gamma=0$, $0<|\gamma|<1$, and $|\gamma|=1$ separately. In each relevant case we determine a necessary condition on the set $\ac(A)$.

Observe $|\gamma| > 1$ yields a contradiction: the expression on the left of \eref{eqn: hyperbola} is nonnegative while the expression on the right of \eref{eqn: hyperbola} is strictly negative. 

If $\gamma = 0$, then the claim is established and it remains to find $\ac(A)$. Then $\lambda_1 = -\lambda_2$ and the system \eref{system: x and y} reduces to $x^2 - y^2 = \frac{1}{2}$. Thus, $|\omega| \in \left[\frac1{\sqrt{2}}, \infty \right)$ and the entries of $B$ all have modulus $\left|(\lambda_2 - \lambda_1)\frac{\omega}2\right| = \rho(A)|\omega|$, i.e., $\ac(A) \subseteq \left[\frac{\rho(A)}{\sqrt{2}},\infty\right)$.

Suppose that $|\gamma| = 1$. Note that $\gamma\not=\pm1$ since this would imply $\lambda_1 = 0$ or $\lambda_2 = 0$. Thus, $\Im(\gamma)\not=0$ and \eref{system: x and y} implies $x = y = 0$, so $\omega = 0$. Hence, the entries of $B$ all have modulus $\left|\frac{\gamma(\lambda_2 - \lambda_1)}{2}\right| = \left|\frac{\lambda_1 + \lambda_2}{2}\right|$, i.e., $\ac(A)\subseteq \left\{ \left|\frac{\lambda_1 + \lambda_2}{2} \right| \right\}$. Since $\gamma^2\not= 1$, it follows that $\Re(\gamma^2) < |\gamma|^4 \leq 1$, as desired.

Suppose that $0 < |\gamma|< 1$. Note that if $\Im(\gamma)=0$, then $x = 0$ by \eref{eqn: line}. Consequently, \eref{eqn: hyperbola} simplifies to $y^2 = -\frac{1 + |\gamma|^2}{2}$, a contradiction. Thus, $\Im(\gamma)\not=0$ and so \eref{eqn: line} is the equation of a line that intersects the origin with slope $-\frac{\Re(\gamma)}{\Im(\gamma)}$. Since $|\gamma|\not=1$, \eref{eqn: hyperbola} is the equation of a hyperbola whose asymptotes have slope $\pm\sqrt{\frac{1 + |\gamma|^2}{1 - |\gamma|^2}}$. A solution to \eref{system: x and y} therefore requires $\left|\frac{\Re(\gamma)}{\Im(\gamma)}\right| < \left|\sqrt{\frac{1 + |\gamma|^2}{1 - |\gamma|^2}}\right|$ or equivalently $0 < \Im(\gamma)^2(1+|\gamma|^2) - \Re(\gamma)^2(1-|\gamma|^2)$. Regrouping gives $0 < \Im(\gamma)^2 - \Re(\gamma)^2 + (\Re(\gamma)^2 + \Im(\gamma)^2)|\gamma|^2 = \Im(\gamma)^2 - \Re(\gamma)^2 + |\gamma|^4$. Utilizing the identity $\Re(\gamma)^2 - \Im(\gamma)^2 = \Re(\gamma^2)$ gives the desired inequality $\Re(\gamma^2)<|\gamma|^4\leq 1$.
To obtain a necessary condition on $\ac(A)$, substitute $y = -\frac{\Re(\gamma)}{\Im(\gamma)} x$ (obtained from \eref{eqn: line}) into \eref{eqn: hyperbola} and solve for $x^2$:
\begin{align*}
x^2 &= \frac{1 - |\gamma|^4}{2\left(1 + |\gamma|^2 - (1 - |\gamma|^2)\frac{\Re(\gamma)^2}{\Im(\gamma)^2} \right)} \\
&= \frac{\Im(\gamma)^2(1 - |\gamma|^4)}{2\left(\Im(\gamma)^2 - \Re(\gamma)^2 + |\gamma|^2(\Im(\gamma)^2 +\Re(\gamma)^2) \right)} \\
&=\frac{\Im(\gamma)^2(1 - |\gamma|^4)}{2(|\gamma|^4 - \Re(\gamma^2))}.
\end{align*}
Thus,
\[
\omega = x + y\ii = \left(1 - \frac{\Re(\gamma)}{\Im(\gamma)} \ii \right)x = \pm\left(\Im(\gamma) - \Re(\gamma)\ii\right) \sqrt{\frac{1 - |\gamma|^4}{2(|\gamma|^4 - \Re(\gamma^2))}} = \pm \gamma\sqrt{\frac{1 - |\gamma|^4}{2(|\gamma|^4 - \Re(\gamma^2))}}\ \ii.
\]
Therefore, the entries of $B$ all have modulus 
\[
\left|(\lambda_2 - \lambda_1)\frac{\gamma - \omega}2\right| = \left|\frac{\lambda_2 - \lambda_1}{2}\gamma\right|\left|1 \pm \sqrt{\frac{1 - |\gamma|^4}{2(|\gamma|^4 - \Re(\gamma^2))}}\ \ii\right| = \left|\frac{\lambda_1 + \lambda_2}{2}\right|\sqrt{1+\frac{1-|\gamma|^4}{2(|\gamma|^4 - \Re(\gamma^2))}},
\]
i.e., $\ac(A) \subseteq \left\{\ \left|\frac{\lambda_1 + \lambda_2}{2}\right|\sqrt{1+\frac{1-|\gamma|^4}{2(|\gamma|^4 - \Re(\gamma^2))}} \right\}$.

We now prove the converse. Assume that $\gamma = 0$ or $\Re(\gamma^2) < |\gamma|^4 \leq 1$. Then $\gamma = 0$, $|\gamma| = 1$ or $0 < |\gamma| < 1$. Let $\kappa\in \left[\frac{\rho(A)}{\sqrt{2}},\infty\right)$ and
\[
\omega = 
\begin{cases}
\sqrt{\frac{1}{2}\big(\frac{\kappa}{\rho(A)}\big)^2 + \frac14} + \sqrt{\frac{1}{2}\big(\frac{\kappa}{\rho(A)}\big)^2 - \frac14}\ \ii & \text{if } |\gamma| = 0\\
0 & \text{if } |\gamma| = 1\\
\gamma\sqrt{\frac{1 - |\gamma|^4}{2(|\gamma|^4 - \Re(\gamma^2))}}\ \ii & \text{if } 0 < |\gamma| < 1.
\end{cases}
\]
Pick a value $b \in \C$ such that $b^2 = \frac{\omega^2 - 1}{4}$. Consider the matrix
$M=\left[\begin{array}{cc}
a & b \\ c & d
\end{array}\right]$,
where $a = 1$, 
$c = \frac{\omega - 1}{2b}$ and $d = \frac{\omega + 1}{2}$. Observe that $\det(M) = 1$ and $\omega = 2bc + 1$. Following the same calculations used to obtain \eref{eqn: mtx B} we find that
\[
B = MAM^{-1} = (\lambda_2-\lambda_1)
\left[\begin{array}{cc}
\frac{\gamma-\omega}{2} & ab \\
-cd & \frac{\gamma + \omega}{2} 
\end{array}\right]
=(\lambda_2-\lambda_1)\left[\begin{array}{cc}
\frac{\gamma-\omega}{2} & b \\
-b & \frac{\gamma + \omega}{2} 
\end{array}\right].
\]
We consider the cases $\gamma = 0$, $|\gamma| = 1$ and $0 < |\gamma| < 1$ separately.

Suppose that $\gamma = 0$. Then $\lambda_2 = -\lambda_1$, $\omega = \sqrt{\frac{1}{2}\big(\frac{\kappa}{\rho(A)}\big)^2 + \frac14} + \sqrt{\frac{1}{2}\big(\frac{\kappa}{\rho(A)}\big)^2 - \frac14}\ \ii$ and
\[
B = 
\lambda_1\left[\begin{array}{cc}
\omega & -2b\\
2b & -\omega
\end{array}\right].
\]
For brevity, write $x = \Re(\omega)$ and $y = \Im(\omega)$. Since $x^2 - y^2 = \frac12$,
\[
|2b|^2 = \left|\omega^2 - 1\right| = |x^2 - y^2 +2xy\ii - 1| = \left|-\frac12 + 2xy\ii \right| = \left|\frac12 + 2xy\ii\right| = |x^2 - y^2 + 2xy\ii| = |\omega|^2.
\]
Therefore, $B$ is uniform with common modulus $|\lambda_1||\omega| = \kappa$, i.e., $A$ is apportionable and $\ac(A) \supseteq \left[\frac{\rho(A)}{\sqrt{2}}, \infty\right)$. It follows that $\ac(A) = \left[\frac{\rho(A)}{\sqrt{2}}, \infty\right)$. 

Suppose that $|\gamma| = 1$. Then $\omega = 0$ and (choosing $b=\frac{\ii}{2}$)
\[
B = 
(\lambda_2-\lambda_1)\left[\begin{array}{cc}
\frac{\gamma}{2} & \frac{\ii}2 \\
-\frac{\ii}2 & \frac{\gamma}{2} 
\end{array}\right] 
= \frac{\lambda_1 + \lambda_2}2\left[\begin{array}{cc}
1 & \frac{\ii}{\gamma} \\
-\frac{\ii}{\gamma} & 1
\end{array}\right].
\]
Therefore, $B$ is uniform with common modulus $\left|\frac{\lambda_1 + \lambda_2}{2}\right|$, i.e., $A$ is apportionable and $\left|\frac{\lambda_1 + \lambda_2}{2}\right| \in \ac(A)$. Since $\sqrt{1+\frac{1-|\gamma|^4}{2(|\gamma|^4 - \Re(\gamma^2))}} = 1$ we have $\ac(A) = \set{\left|\frac{\lambda_1 + \lambda_2}{2}\right|\sqrt{1+\frac{1-|\gamma|^4}{2(|\gamma|^4 - \Re(\gamma^2))}}}$.

Finally, suppose that $0 < |\gamma| < 1$. Then $\omega = \gamma\sqrt{\frac{1 - |\gamma|^4}{2(|\gamma|^4 - \Re(\gamma^2))}}\ \ii$ and 
\[
B =
\frac{\lambda_2-\lambda_1}2\left[\begin{array}{cc}
\gamma-\omega & 2b \\
-2b & \gamma + \omega 
\end{array}\right].
\]
It is not difficult to see that the $(1,2)$-entry and $(2,1)$-entry of $B$ have the same modulus. The $(1,1)$-entry and $(2,2)$-entry of $B$ have the same modulus since
\[
|\gamma - \omega| = |\gamma|\left|1 - \sqrt{\frac{1 - |\gamma|^4}{2(|\gamma|^4 - \Re(\gamma^2))}}\ \ii\right| = |\gamma|\left|1 + \sqrt{\frac{1 - |\gamma|^4}{2(|\gamma|^4 - \Re(\gamma^2))}}\ \ii\right| = |\gamma + \omega|.
\]
To establish the $(1,1)$-entry and $(1,2)$-entry have the same modulus, we first note that
\[
(\gamma\overline{\omega})^2 = \left(\gamma\overline{\gamma}\sqrt{\frac{1 - |\gamma|^4}{2(|\gamma|^4 - \Re(\gamma^2))}}\ (-\ii)\right)^2 = -|\gamma|^4\frac{1 - |\gamma|^4}{2(|\gamma|^4 - \Re(\gamma^2))}.
\]
Then by Observation~\ref{obs: 2x2 distinct evals},
\begin{align*}
\frac{|1 - \omega^2|^2 - |\gamma^2 - \omega^2|^2}2 &= \frac{1 - |\gamma|^4}{2} - \Re(\omega^2) + \Re((\gamma\overline{\omega})^2)\\
&=\frac{1 - |\gamma|^4}{2} + \Re(\gamma^2) \frac{1 - |\gamma|^4}{2(|\gamma|^4-\Re(\gamma^2))} - |\gamma|^4 \frac{1 - |\gamma|^4}{2(|\gamma|^4-\Re(\gamma^2))}\\
&=\frac{1 - |\gamma|^4}{2}\left(1 + \frac{\Re(\gamma^2) - |\gamma|^4}{|\gamma|^4-\Re(\gamma^2)} \right)\\
&= 0.
\end{align*}
So, $|1 - \omega^2| = |\gamma^2 - \omega^2|$. It follows that 
\[ |2b|^2 = |\omega^2-1| = |\gamma^2 - \omega^2| = |\gamma-\omega||\gamma+\omega| = |\gamma-\omega|^2.
\]
Thus, $B$ is uniform with common modulus
\[
\left|\frac{\lambda_2-\lambda_1}2 \right| \left|\gamma - \omega \right| = \left|\frac{\lambda_2-\lambda_1}2 \right| \left|\frac{\lambda_2 + \lambda_1}{\lambda_2 - \lambda_1}\right| \left|1 - \sqrt{\frac{1 - |\gamma|^4}{2(|\gamma|^4 - \Re(\gamma^2))}}\ \ii \right| = \left|\frac{\lambda_1 + \lambda_2}2 \right| \sqrt{1 + \frac{1 - |\gamma|^4}{2(|\gamma|^4 - \Re(\gamma^2))}},
\]
i.e., $A$ is apportionable and $\left|\frac{\lambda_1 + \lambda_2}2 \right| \sqrt{1 + \frac{1 - |\gamma|^4}{2(|\gamma|^4 - \Re(\gamma^2))}} \in \ac(A)$. Since this gives containment in both directions, $\ac(A) = \left\{\left|\frac{\lambda_1 + \lambda_2}2 \right| \sqrt{1 + \frac{1 - |\gamma|^4}{2(|\gamma|^4 - \Re(\gamma^2))}}\right\}$. 
\end{proof}

\begin{remark}
Theorem~\ref{thm: 2x2 distinct evals} demonstrates that the bound in Lemma~\ref{lem: Had lower bound} is sharp for some matrices. Indeed, if $A\in\C^{2\times 2}$ with $\gamma = 0$, then $\ac(A) = \left[\frac{\rho(A)}{\sqrt{2}}, \infty\right) = \left[\frac{|\det(A)|^{1/2}}{\sqrt{2}}, \infty\right)$. Interestingly, Lemma~\ref{lem: GL lower bound} is not sharp in this case.
\end{remark}

The following corollary builds upon Theorem~\ref{thm: 2x2 distinct evals} by providing easily verifiable necessary and sufficient conditions that a nonsingular $2\times 2$ matrix is apportionable.

\begin{corollary}\label{cor: 2x2 nonsing}
Let $A\in\C^{2\times 2}$ have nonzero eigenvalues $\lambda_1 = r_1 e^{\ii\theta_1}$ and $\lambda_2 = r_2 e^{\ii\theta_2}$, where $r_i > 0$ and $0\leq \theta_i < 2\pi$ for $i\in\{1,2\}$. Let $r = \frac{r_2}{r_1}$ and $\theta = \theta_2 - \theta_1$.
Then $A$ is apportionable if and only if $\lambda_1 = -\lambda_2$, $\lambda_1 = c\ii \lambda_2$ for some $c\in\R$, or $\frac{\pi}{2} < \theta < \frac{3\pi}{2}$ and $|r\cos(\theta) + 1| < |\sin(\theta)|$.
\end{corollary}
\begin{proof}
Assume that $A$ is apportionable (hence $\lambda_1\not=\lambda_2$) and let $\gamma = \frac{\lambda_2 + \lambda_1}{\lambda_2 - \lambda_1}$. Then $\gamma = \frac{\mu + 1}{\mu - 1}$, where $\mu = \frac{\lambda_2}{\lambda_1} = re^{\ii \theta}$. By Theorem~\ref{thm: 2x2 distinct evals} $\mu = -1$ or $\Re(\gamma^2) < |\gamma|^4 \leq 1$. In the case $\mu = -1$, by definition $\lambda_1 = -\lambda_2$. Now consider the case $\mu\not= -1$ so that $\Re(\gamma^2) < |\gamma|^4 \leq 1$. Let $x = r\cos(\theta)$ and $y = r\sin(\theta)$ so that $\mu = x + y\ii$. Direct computation yields
\[
\Re(\gamma^2) 
= \Re\left(\left(\frac{(\mu+1)(\overline{\mu}-1)}{|\mu-1|^2}\right)^2\right)
=\Re\left(\frac{(|\mu|^2-1-2\ii y)^2}{|\mu-1|^4}\right)
=\frac{(|\mu|^2-1)^2-4y}{|\mu-1|^4}.
\]
Since $|\gamma|^4 = \frac{|\mu + 1|^4}{|\mu - 1|^4}$, the inequality $\Re(\gamma^2) < |\gamma|^4 \leq 1$ is equivalent to 
\begin{equation}\label{eqn: cor 2x2}
(|\mu|^2-1)^2-4y^2<|\mu+1|^4 \leq |\mu - 1|^4.
\end{equation}
The rightmost inequality of \eref{eqn: cor 2x2} implies $x \leq 0$. If $x = 0$, then the definition of $\mu$ implies $\lambda_1 = c\ii\lambda_2$ for some $c\in\R$. So, suppose $x = r\cos(\theta) < 0$ which implies $\frac{\pi}{2} < \theta < \frac{3\pi}{2}$. Simplifying the leftmost inequality of \eref{eqn: cor 2x2} using the identity $x^2 + y^2 = r^2$ yields 
\begin{equation}\label{eqn: cor 2x2 second}
0 < 2r^2 + r\cos(\theta)(r^2 + 1).
\end{equation}
Multiply both sides by $\frac{\cos(\theta)}{r} < 0$ and simplify to obtain $0 > (r\cos(\theta) + 1)^2 - \sin^2(\theta)$. Thus, $|\sin(\theta)| > |r\cos(\theta) + 1|$, as desired.

The converse follows a similar argument. We note that $\frac{\pi}{2} < \theta < \frac{3\pi}{2}$ implies $\lambda_1\not=\lambda_2$.
\end{proof}

\begin{figure}[H]
\begin{center}
\includegraphics[scale=0.6]{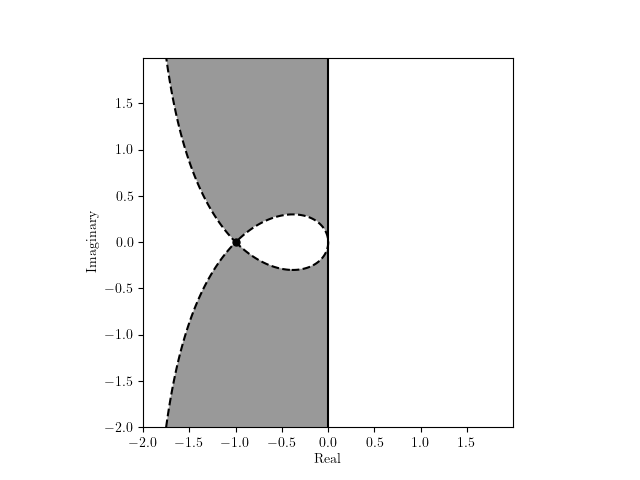}
\end{center}
\caption{Admissible region for the second eigenvalue of a 2×2 apportionable matrix with one eigenvalue equal to 1.
}\label{fig:evals2by2}
\end{figure}

Figure~\ref{fig:evals2by2} provides a visualization of the admissible eigenvalues of a $2\times 2$ apportionable matrix; we fix an eigenvalue to be $1$ and plot the region that contains the remaining eigenvalue. See \cite{code} for a Python script that produces an interactive image that allows the user to select the value of an eigenvalue. 

We now use Theorem~\ref{thm: 2x2 distinct evals} and Corollary~\ref{cor: 2x2 nonsing} to address \cite[Question 9.2]{CCGH2024} which asks whether the set of apportionable (or non-apportionable) matrices has a positive Lebesgue measure. While not fully explored in this paper, we have observed that for $n=2$, there exists a set of apportionable matrices and a set of non-apportionable matrices that have positive Lebesgue measure. The next example shows that all $2 \times 2$ matrices sufficiently close to the identity are not apportionable.

\begin{example}\label{ex: non-apport pos meas}
Let $B\in \C^{2\times 2}$ such that $\|I - B\|_F < \frac12$. Let the eigenvalues of $B$ be denoted by $\lambda_1$ and $\lambda_2$. The eigenvalues of $I - B$ are $1 - \lambda_1$ and $1 - \lambda_2$. If $\lambda_1 = \lambda_2$, then $I - B$ is not apportionable, so suppose $\lambda_1\not= \lambda_2$. Since the Frobenius norm is unitarily invariant, and $I-B$ is unitarily similar to an upper triangular matrix, $\|I - B\|_F = \sqrt{|1-\lambda_1|^2 + |1-\lambda_2|^2 + |a|^2} < \frac12$, where $a\in\C$ is the $(1,2)$-entry of the upper triangular matrix. It follows that $|1-\lambda_j| < \frac{1}{2}$ for $j\in\{1,2\}$, which implies $\Re(\lambda_j) > \frac12$. Thus, $|\lambda_2 + \lambda_1| > 1$. By the triangle inequality, we also have $|\lambda_2 - \lambda_1| \leq |\lambda_2 - 1| + |1 - \lambda_1| < 1$. Let $\gamma = \frac{\lambda_2 + \lambda_1}{\lambda_2 - \lambda_1}$. Since $|\gamma| > 1$, Theorem~\ref{thm: 2x2 distinct evals} guarantees $B$ is not apportionable. Therefore, the open ball of radius $\frac12$ centered at the identity matrix in $\C^{2\times 2}$ provides an example of a set of non-apportionable matrices with positive Lebesgue measure.
\end{example}

We now use the Gershgorin Circle Theorem to obtain a set of apportionable $2\times2$ matrices with positive Lebesgue measure. Recall, for any $A = [a_{i,j}]\in\C^{n\times n}$ the Gershgorin Circle Theorem states that every eigenvalue $\lambda$ of $A$ lies in one of the discs $\{z \in \C : |z - a_{i,i}| \leq \sum_{j\not=i}|a_{i,j}|\}$ for $i\in \{1,\ldots,n\}$. Furthermore, if these discs are disjoint, then each disc contains exactly one eigenvalue.

\begin{example}
Let $A = \diag\big(1,-\frac12 + \ii\big)$ and $B = [b_{i,j}]$ be a $2\times 2$ matrix such that $\|A - B\|_F < \frac{1}{12}$. Let the eigenvalues of $B$ be denoted by $\lambda_1$ and $\lambda_2$. The Gershgorin Circle Theorem implies that, without loss of generality, $|\lambda_1 - b_{1,1}| < |b_{1,2}|$ and $|\lambda_2 - b_{2,2}| < |b_{2,1}|$. Thus, $|\lambda_1 - 1| - |1 - b_{1,1}| \leq |\lambda_1 - b_{1,1}| < \frac1{12}$ and so $|\lambda_1 - 1| < \frac16$. Similarly, $|\lambda_2 + \frac12 - \ii| < \frac16$. 

As in the proof of Corollary~\ref{cor: 2x2 nonsing}, let $\mu = \frac{\lambda_2}{\lambda_1} = x + y\ii$, where $x = r\cos(\theta)$ and $y = r\sin(\theta)$. For the case $x < 0$ in the proof of Corollary~\ref{cor: 2x2 nonsing} we obtain \eref{eqn: cor 2x2 second} which we restate:  $0 < 2r^2 + r\cos(\theta)(r^2 + 1)$. Expressing \eref{eqn: cor 2x2 second} in terms of $x$ and $y$ yields
\[
x < 0 < (x+1)\parens{\parens{x+\frac12}^2 +y^2 -\frac14} + y^2.
\] 
This gives rise to a sufficient condition for apportionability: 
\[
-1 < x < 0\quad\text{and}\quad |y|>\frac12.
\] 
Hence, it suffices to show $\big|\mu - \big(-\frac12 + \ii \big)\big| < \frac12$. Indeed
\begin{align*}
\Big|\frac{\lambda_2}{\lambda_1} - \parens{-\frac12 + \ii}\Big| &= \frac{\left|\left(\lambda_2-\left(-\frac12 + \ii\right)\right) - \left(\lambda_1 - 1\right)\left(-\frac12 + \ii\right)\right|}{|\lambda_1|} \\
&\leq \frac{\left|\lambda_2-\left(-\frac12+\ii\right)\right| + \left|\lambda_1 - 1\right|\frac{\sqrt{5}}{2}}{|\lambda_1|} \\
&< \frac{\frac16\left(1 + \frac{\sqrt{5}}{2}\right)}{\frac56} \\
&< \frac12.
\end{align*}
Therefore, the open ball of radius $\frac1{12}$ centered at $A$ in $\C^{2\times 2}$ provides an example of a set of apportionable matrices with positive Lebesgue measure.
\end{example}

\subsection{Apportionable Matrices of order 3}\label{sec: order 3}

This section initiates the classification of apportionable $3\times 3$ matrices.  A summary of all results for $3\times 3$ apportionable matrices $A$ for which $\ac(A)$ is known is presented in Table~\ref{tab: 3 by 3}. Example~\ref{ex: 3x3 lambda I oplus O} and Example~\ref{ex: 3x3 I oplus nilpotent} establish that matrices with Jordan form
\[
\left[\begin{array}{ccc} 
\lambda & 1 & 0 \\
0 & \lambda & 0 \\
0 & 0 & 0 
\end{array}\right]
\qquad \text{and} \qquad
\left[\begin{array}{ccc} 
\lambda & 0 & 0 \\
0 & 0 & 1 \\
0 & 0 & 0 
\end{array}\right],
\]
where $\lambda$ is nonzero, are apportionable. Theorem~\ref{thm: 2x2 distinct evals} and Lemma~\ref{lem:pad-by-0} establish that some matrices with Jordan form
\[
\left[ \renewcommand{\arraystretch}{1} \begin{array}{ccc} 
\lambda_1 & 0 & 0 \\
0 & \lambda_2 & 0 \\
0 & 0 & 0 
\end{array} \right],
\]
where $\lambda_1\not=\lambda_2$ are nonzero, are apportionable. There is nothing known for the remaining cases, listed below ($\lambda_1, \lambda_2$ and $\lambda_3$ are distinct and nonzero):
\[
\left[\begin{array}{ccc} 
\lambda_1 & 1 & 0 \\
0 & \lambda_1 & 1 \\
0 & 0 & \lambda_1 
\end{array} \right], 
\quad
\left[\begin{array}{ccc} 
\lambda_1 & 1 & 0 \\
0 & \lambda_1 & 0 \\
0 & 0 & \lambda_2 
\end{array} \right]
\quad \text{and} \quad
\left[\begin{array}{ccc} 
\lambda_1 & 0 & 0 \\
0 & \lambda_2 & 0 \\
0 & 0 & \lambda_3 
\end{array} \right].
\]
\begin{table}[H]
\begin{center}
\resizebox{0.98\textwidth}{!}{
\renewcommand{\arraystretch}{1.8}
\begin{tabular}{|c||*5{>{\renewcommand{\arraystretch}{1}}c|}}
\hline
Matrix $A$ 
& 
$\left[\begin{array}{ccc} 
0 & 0 & 0 \\
0 & 0 & 0 \\
0 & 0 & 0
\end{array} \right]$ 
& 
$\left[\begin{array}{ccc} 
\lambda_1 & 0 & 0 \\
0 & \lambda_1 & 0 \\
0 & 0 & \lambda_1 
\end{array} \right]$
& 
$\left[\begin{array}{ccc} 
0 & 1 & 0 \\
0 & 0 & 0 \\
0 & 0 & 0 
\end{array} \right]$ 
& 
$\left[\begin{array}{ccc} 
0 & 1 & 0 \\
0 & 0 & 1 \\
0 & 0 & 0
\end{array} \right]$ 
\\[3pt]
\hline
Apportionable 
& Yes & No & Yes & Yes
\\[3pt]
\hline
$\ac(A)$ 
& $\{0\}$ & $\emptyset$ & $(0,\infty)$ & $(0,\infty)$
\\[3pt]
\hline
Proof 
& $A = O$ & $A = \lambda_1 I$ & Thm.~\ref{thm: nilpotent} &  Thm.~\ref{thm: nilpotent}
\\[3pt]
\hline
\hline
Matrix $A$ 
&
$\left[\begin{array}{ccc} 
\lambda_1 & 0 & 0 \\
0 & 0 & 0 \\
0 & 0 & 0
\end{array} \right]$
&
$\left[\begin{array}{ccc} 
\lambda_1 & 1 & 0 \\
0 & \lambda_1 & 0 \\
0 & 0 & \lambda_1 
\end{array} \right]$
&
$\left[\begin{array}{ccc} 
\lambda_1 & 0 & 0 \\
0 & \lambda_1 & 0 \\
0 & 0 & \lambda_2 
\end{array} \right]$ & $\left[\begin{array}{ccc} 
\lambda_1 & 0 & 0 \\
0 & \lambda_1 & 0 \\
0 & 0 & 0 
\end{array} \right]$ 
\\[3pt]
\hline
Apportionable 
& Yes & No & Yes $\Leftrightarrow$ $\Re\big(\tfrac{\lambda_2}{\lambda_1}\big) = -\frac12$ & No 
\\[3pt]
\hline
$\ac(A)$ 
& $\Big[\frac{|\lambda_1|}3, \infty \Big)$ & $\emptyset$
& $\left\{|\lambda_1|\sqrt{\frac{\Im\left(\lambda_2/\lambda_1\right)^2}{(3-2s)^2} +\frac14} : s = 0,1 \right\}$ & $\emptyset$
\\[3pt]
\hline
Proof 
& Prop.~\ref{prop: rank 1 apport} & Thm.~\ref{thm: I oplus lambda} & Thm.~\ref{thm: I oplus lambda} & Thm.~\ref{thm: I oplus lambda}
\\[3pt]
\hline
\end{tabular}
}
\end{center}
\caption{$3\times 3$ apportionment results in Jordan form ($\lambda_1$ and $\lambda_2$ are distinct and nonzero). }
\label{tab: 3 by 3}
\end{table}

We now provide two examples of apportionable $3\times 3$ matrices that are not covered by preceding results. Note that the next example provides an example of a nonsingular matrix $B\in\C^{n\times n}$ for which $0<\varsigma(B)<n$ (c.f. Problem~\ref{prob: add zeros}). Indeed, the upper $2\times 2$ block of $A$ in the next example is $J_2(\lambda)$ which is not apportionable for $\lambda \neq 0$ by \cite[Theorem 7.4]{CCGH2024}.
\begin{example}\label{ex: 3x3 lambda I oplus O}
Let $\lambda \in \C$. We show that the matrix
\[
A = 
\begin{bmatrix}
\lambda & 1 & 0 \\
0 & \lambda & 0 \\
0 & 0 & 0 
\end{bmatrix}
\]
is apportionable. The case $\lambda = 0$ is handled by Theorem \ref{thm: nilpotent}, so we may assume $\lambda \neq 0$. Consider the matrices
\[
M = 
\left[\begin{array}{ccc}
0 & 1 & 1 \\
e^{2\pi \ii/3} & 0 & 1 \\
1 & 0 & 0
\end{array}\right]
\left[\begin{array}{ccc}
\lambda & 0 & 0 \\
0 & 1 & 0 \\
0 & 0 & 1 
\end{array}\right]
\qquad \text{and} \qquad
M^{-1} =
\left[\begin{array}{ccc}
\lambda^{-1} & 0 & 0 \\
0 & 1 & 0 \\
0 & 0 & 1 
\end{array}\right]
\left[\begin{array}{ccc}
0 & 0 & 1 \\
1 & -1 & e^{2\pi \ii/3} \\
0 & 1 & -e^{2\pi \ii/3} 
\end{array}\right].
\]
Straightforward computation yields 
\[
MAM^{-1} = 
\left[\begin{array}{ccc}
0 & 1 & 1 \\
e^{2\pi \ii/3} & 0 & 1 \\
1 & 0 & 0 \\
\end{array}\right]
\left[\begin{array}{ccc}
\lambda & \lambda & 0 \\
0 & \lambda & 0 \\
0 & 0 & 0 
\end{array}\right]
\left[\begin{array}{ccc}
0 & 0 & 1 \\
1 & -1 & e^{2\pi \ii/3} \\
0 & 1 & -e^{2\pi \ii/3} 
\end{array}\right] 
= \lambda\left[\begin{array}{ccc}
1 & -1 & e^{2\pi \ii/3} \\
e^{2\pi \ii/3} & -e^{2\pi \ii/3} & -1 \\
1 & -1 & 1+e^{2\pi \ii/3} 
\end{array}\right].
\]
Then $MAM^{-1}$ is uniform with a common modulus $|\lambda|$. Thus, $A$ is apportionable.

\end{example}

\begin{example}\label{ex: 3x3 I oplus nilpotent}
 Let $\lambda \in \C$. We show that the matrix
\[
A = 
\begin{bmatrix}
    \lambda & 0 & 0 \\
    0 & 0 & 1 \\
    0 & 0 & 0 
\end{bmatrix}
\]
is apportionable. The case $\lambda = 0$ is handled by Theorem \ref{thm: nilpotent}, so we may assume $\lambda \neq 0$. Consider the matrices
\[
M = 
\left[\begin{array}{ccc}
    0 & 1 & e^{2\pi \ii/3} \\
    1 & 0 & e^{\pi \ii/3} \\
    1 & 1 & 0
\end{array}\right]
\left[\begin{array}{ccc}
    1 & 0 & 0 \\
    0 & \lambda & 0 \\
    0 & 0 & 1 
\end{array}\right]
\text{ and } 
M^{-1} =
\frac{1}{1 + e^{\pi\ii/3}}
\left[\begin{array}{ccc}
    1 & 0 & 0 \\
    0 & \lambda^{-1} & 0 \\
    0 & 0 & 1 
\end{array}\right]
\left[\begin{array}{ccc}
    -1 & e^{\pi\ii/3} & 1 \\
    1 & -e^{\pi\ii/3} & e^{\pi\ii/3} \\
    e^{-\pi\ii/3} & e^{-\pi\ii/3} & -e^{-\pi\ii/3}
\end{array}\right].
\]
Straightforward computation yields

\begin{align*}
MAM^{-1} & =
\frac{1}{1 + e^{\pi\ii/3}}
\left[\begin{array}{ccc}
    0 & 1 & e^{2\pi\ii/3} \\
    1 & 0 & e^{\pi \ii/3} \\
    1 & 1 & 0
\end{array}\right]
\left[\begin{array}{ccc}
    \lambda & 0 & 0 \\
    0 & 0 & \lambda \\
    0 & 0 & 0 
\end{array}\right]
\left[\begin{array}{ccc}
    -1 & e^{\pi\ii/3} & 1 \\
    1 & -e^{\pi\ii/3} & e^{\pi\ii/3} \\
    e^{-\pi\ii/3} & e^{-\pi\ii/3} & -e^{-\pi\ii/3}
\end{array}\right] \\
& = \frac{\lambda}{1 + e^{\pi\ii/3}}\left[\begin{array}{ccc}
    1 & 1 & -1 \\
    -e^{\pi\ii/3} & e^{2\pi\ii/3} & e^{\pi\ii/3} \\
    1-e^{\pi\ii/3} & 1+e^{2\pi\ii/3} & e^{\pi\ii/3} - 1 
\end{array}\right].
\end{align*}
Then $MAM^{-1}$ is uniform with a common modulus $\frac{1}{\sqrt{3}}|\lambda|$. Thus, $A$ is apportionable.
\end{example}

\section{Conclusion}\label{sec: conclusion}

This paper complements the initial exploration in \cite{CCGH2024}, which primarily focused on $\U$-apportionable matrices, by shifting the focus to the general case. Our primary contributions include establishing the apportionability of all nilpotent matrices and proving that any $n \times n$ matrix $A \in \mathbb{C}^{n \times n}$ with $\rank(A) \leq \frac{n}{2}$ is apportionable. We also initiated the study of sets of apportionment constants showing that $\ac(A) = (0,\infty)$ for any nonzero nilpotent matrix $A$ and identifying nonzero matrices with finite sets of apportionment constants.

The characterization of apportionable matrices of small order provided critical insights and yielded many of the techniques and tools developed in this paper. We achieved a complete classification for $2 \times 2$ matrices and made substantial progress towards a classification for $3 \times 3$ matrices. Further efforts to complete the characterization of all $3 \times 3$ apportionable matrices would likely uncover additional general results. Based on our study of $2 \times 2$ matrices, it appears that the case involving three distinct eigenvalues will present the most significant challenge in the $3 \times 3$ classification.

Several open problems and avenues for future research emerge from this work. We are particularly interested in the relationship between the rank and dimension of an apportionable matrix, especially regarding Problem~\ref{prob: add zeros}. A promising starting point for investigating these problems is to determine $\varsigma(I_n)$. For instance, it remains unknown whether $\varsigma(I_3)$ is $2$ or $3$. 

Although Section 5.1 has fully characterized apportionability of $2\times 2$ matrices, another avenue of future research is characterizing $\varsigma(A)$ for all matrices $A\in\C^{2\times 2}$. 
We found that the set of apportionable matrices and the set of non-apportionable both have positive measure. This naturally leads to the question: When we partition the non-apportionable matrices as $\{A \in \C^{2\times2} : \varsigma(A) = 1\}\cup\{A \in \C^{2\times2} : \varsigma(A) = 2\}$, does each set have positive measure? Note that the white region in Figure \ref{fig:evals2by2} consists of all $\lambda \in \C$ such that the matrix $A = \diag(1,\lambda)$ has either $\varsigma(A) = 1$ or $\varsigma(A) = 2$.

A question related to the study of $\varsigma$ is whether Theorem \ref{thm: half rank apportionable} is sharp. That is, what is the least $r > \frac{n}{2}$ such that there exists a non-apportionable $A \in \mathbb{C}^{n \times n}$ with $\rank(A)=r$? For $3 \times 3$ matrices, the answer is $r=2$ as shown by $I_2 \oplus O_1$ which is not apportionable (see Table \ref{tab: 3 by 3}). This suggests that simple candidates for such non-apportionable matrices might be $I_r \oplus O_{n-r}$ for any order $n$ and rank $r > \frac{n}{2}$.

Beyond the classification of matrices and their apportionment constants, another avenue for research involves studying the set of matrices that apportion a fixed matrix. Understanding the nature of these transforming matrices may provide further insight into apportionment.

\section{Acknowledgements}
The research of D.\ Baker, B.\ Curtis, J.\ Miller and H.\ Pungello was partially supported by NSF grant 1839918 and the authors thank the National Science Foundation.

This paper was approved for public release: AFRL-2025-4326. The views expressed in this article are those of the authors and do not necessarily reflect the official policy or position of the Air Force, the Department of Defense, or the U.S. Government.


\end{document}